\documentclass[10pt]{article}
\usepackage{amssymb}
\usepackage{amsmath}

\begin{document}
\begin{center}
\Large
% TITLE GOES HERE
Parametric restrictions on quasi-symmetric designs
\end{center}

\begin{flushleft}
Bhaskar Bagchi\footnote{The author is a retired professor of Indian Statistical Institute.}  \\
1363, 10th Cross Road,\\
Kengeri Satellite Town, \\
Bangalore 560 060, \\
India.\\
\verb+bhaskarbagchi53@gmail.com+
\end{flushleft}

\begin{abstract}
In this paper, we attach several new invariants to connected strongly regular graphs (excepting conference graphs on non-square number of vertices) : one invariant called the discriminant, and a p-adic invariant corresponding to each prime number p. We prove parametric restrictions on quasi-symmetric 2-designs with a given connected block graph $G$ and a given defect (absolute difference of the two intersection numbers) solely in terms of the defect and the parameters of $G$, including these new invariants.  This is a natural analogue of Schutzenberger's Theorem and the Shrikhande-Chowla-Ryser theorem. This theorem is effective when these graph invariants can be explicitly computed. We do this for complete multipartite graphs, co-triangular graphs, symplectic non-orthogonality graphs (over the field of order $2$) and the  Steiner graphs, yielding explicit restrictions on the parameters of quasi-symmetric 2-designs whose block graphs belong to any of these four classes.
\end{abstract}
\newline

\newpage

\section{\textbf{Introduction}}
An \textbf{ incidence system} is a pair $(P, \mathbb{B})$ where $P$ is a set whose elements are called the points, and $\mathbb{B}$ is a collection of subsets of $P$, called blocks.

Recall that a \textbf{2-design} with parameters $v,k,\lambda$ (in short a $2-(v,k,\lambda)$ design) is an incidence system with $v$ points in total and with $k$ points in each block such that any two distinct points are together in exactly $\lambda$ blocks. We shall sometimes refer to the parameter $\lambda$ as the \textbf{balance parameter} of the design. An easy counting argument shows that a 2-design has ancillary parameters $b,r$ such that the number of blocks is $b$ and each point is in exactly $r$ blocks. These ancillary parameters are determined by the main parameters by the formulae
\begin{equation} \label{feasible} bk=rv,\; r(k-1) = \lambda (v-1). \end{equation}
We say that $b,v,r,k, \lambda$ are \textbf{feasible parameters}  for a 2-design if these are positive integers satisfying (\ref{feasible}) (even if the design may not exist). A 2-design is non-trivial if $v >k$, equivalently, if $r > \lambda$. The number $r-\lambda$ is called the \textbf{order} of the 2-design.

 A well known theorem of Fisher says that the parameters of a non-trivial 2-design satisfy $b \geq v$, equivalently $r \geq k$. A 2-design with $b=v$, $r=k$ is called a \textbf{symmetric 2-design}. Evidently, a $2-(v,k, \lambda)$ design is symmetric if and only if $\lambda (v-1)=k(k-1)$. The statistician's terminology for a 2-design is `balanced incomplete block design' (BIBD) and a symmetric 2-design is called a `symmetrical balanced incomplete block design' (SBIBD) in the Statistics literature. A good source on Design Theory is the monograph \cite{Cameron}.

An \textbf{intersection number} in a design is a number $x$ such that some pair of distinct blocks  have exactly $x$ points in common.
A 2-design is symmetric if and only if it has exactly one intersection number (in which case the intersection number necessarily equals the balance parameter). In view of this result, it is natural to define a \textbf{quasi-symmetric} 2-design to be a design with exactly two intersection numbers. Thus, according to the definition adopted here, the symmetric 2-designs are not quasi-symmetric. A fairly comprehensive source on quasi-symmetric designs is the book \cite{MSS}.

 We denote the intersection numbers of a quasi-symmetric 2-design by $\lambda_1 < \lambda_2$.  The \textbf{defect} of such a design is the difference $\lambda_2 - \lambda_1$. The \textbf{block graph} of a quasi-symmetric 2-design is the graph whose vertices are the blocks of the design, two blocks are adjacent if and only if they have $\lambda_2$ points in common. Some authors define adjacency in the block graph in terms of the smaller intersection number. This is not a serious difference since it only replaces the block graph by its complement. The definition adopted here is consonant with the usual definition of the line graphs of partial linear spaces and graphs. It is also the definition used in \cite{Cameron} (but not in \cite{MSS}).

The graphs considered here are simple graphs. That is, they are loopless, undirected and without multiple edges. A graph is said to be \textbf{regular} of degree $a$ if each vertex has exactly $a$ neighbours. A regular graph is said to be \textbf{strongly regular} (s.r.g.) if there are constants $c,d$ such that any two distinct vertices have exactly $c$ or $d$ common neighbours, according as these two vertices are themselves adjacent or not. A famous theorem of Goethals and Seidel (\cite{Goethals}) says that the block graph of any quasi-symmetric 2-design is a  strongly regular graph. It is easy to see that the only quasi-symmetric 2-designs with disconnected block graphs are the multiples of symmetric 2-designs, obtained by repeating each block a constant number of times. We exclude these designs from consideration. Thus, for us, all the block graphs are connected.

For any graph $G$, its complement  is defined as the graph $G^*$ such that (i) $G^*$ has the same vertices as $G$, and (ii) two distinct vertices are adjacent in $G^*$ if and only if they are non-adjacent in $G$. If $G$ is strongly regular, then so is $G^*$.

The methods used in the next section (Section 2)  of this paper are elementary. The  organization (in terms of the spectrum of a potential block graph and the defect of a putative quasi-symmetric 2-design) of the results here may have its uses in furthering the subject.
 Recall that the complement $\overline{D}$ of an incidence system $D$ is the incidence system whose blocks are the complements (relative to the point set) of the blocks of $D$. It is easy to see that the complement of a 2-design is a 2-design with the same order. Further, the complement of a  quasi-symmetric 2-design with block graph $G$ and defect $\mu$ is again a quasi-symmetric 2-design with block graph $G$ and defect $\mu$. (More precisely, the two block graphs are isomorphic, and an isomorphism between them is given by complementation.) In the first result (Theorem 2.1) of Section 2, we observe that, up to complementation, the parameters of a quasi-symmetric 2-design (including the two intersection numbers) are determined by the defect of the design and the parameters of the block graph. Therefore, it is natural to formulate the parametric restrictions on quasi-symmetric 2-designs in terms of the graph parameters and the defect. We point out that it is not true (as has sometimes been stated) that the graph parameters alone determine the design parameters. Even aside from the ambiguity due to complementation, a famous construction due to Shrikhande and Raghavarao \cite{SSS4} provides non-trivial counter examples. See Theorem 4.4 below for a convenient reformulation of the  result of \cite{SSS4}.

 As an easy  consequence of Theorem 2.1, we present  Theorem 2.3, which gives elementary parametric restrictions on quasi-symmetric 2-designs. Given an  s.r.g. $G$ and a positive integer $\mu$,  Theorem 2.3 gives necessary and sufficient conditions on the parameters of   $G$ for the design parameters given by Theorem 2.1 to be feasible  for a quasi-symmetric 2-design of defect $\mu$. (This, of course, does not guarantee that the design exists! See Definition 2.2 for the precise meaning of feasibility in this context.) The  corollaries to Theorem 2.3 assert the finiteness of certain families of quasi-symmetric 2-designs of a given defect $\mu \geq 2$.
 Corollary 2.4 is a complete classification of the feasible parameters of  quasi-symmetric $2$-designs with complete multipartite block graphs. This result may look complicated, but it lends itself to a fast computation of tables of feasible parameters of such designs, even by hand. For instance, this result has the surprising consequence that, given $n \geq 2, \mu \geq 1$, there is at most one value of $m\geq 2$ for which a quasi-symmetric 2-design of defect $\mu$ and block graph $K_{m \times n}$ is feasible. Namely, if $n=2$ then $m=4\mu-1$. If $n >2$ then, writing $\mu$ (uniquely) as $\mu = (n-1)t + \alpha$ ($t \geq 0$, $0 \leq \alpha < n-1$), we get that if $\alpha >0$, $(n-2\alpha)^2 -4 \alpha (\alpha -1)$ is a perfect square and $\alpha$ divides $n^2t$ then $m= \frac{n^2t}{\alpha} + n+1$ (otherwise there is no such $m$).

  In the context of Section 2, we should recall that, in \cite{Pawale1}, Pawale classified the feasible parameters of quasi-symmetric 2-designs of defect $\mu = 1$. Up to complementation, these are the Steiner $2$-designs, and the residuals of biplanes. In \cite{Pawale2} and \cite{Mavron}, the feasible parameters for q.s. $2$-designs with $\mu = 2,3$ are classified. The recent non-existence results in \cite{Pawale3} and \cite{Pawale4} are consequences of feasibility; as such, they do not rule out any feasible parameters. Likewise, Section 2 handles the question of feasibility only. In contrast, the main result in Section 4 rules out infinitely many feasible parameters.

A central concern of design theory is the following hard question : given parameters $v,k,\lambda$ satisfying the necessary condition $k(k-1) = \lambda (v-1)$, when does a (symmetric) $2-(v,k,\lambda)$ design exist? In \cite{Schutzen}, Schutzenberger proved :
\newline

\textbf{Theorem 1.1}

 If $v$ is an even number then, for the existence of a symmetric 2-design of order $\nu$ on $v$ points, $\nu$ must be a perfect square.

In \cite{SSS1} and \cite{Chowla}, Shrikhande and Chowla-Ryser generalized a previous result of Bruck and Ryser (\cite{Bruck}, the case $\lambda=1$) to prove :
\newline

\textbf{Theorem 1.2}

 If $v$ is an odd number, then for the existence of a symmetric 2-design on $v$ points with balance $\lambda$ and order $\nu$, the diophantine equation $\nu x^2 + (-1)^{(v-1)/2} \lambda y^2 = z^2$ must have a solution $(x,y,z)\not= (0,0,0)$ in integers. \newline

Apart from an isolated computer-assisted non-existence result (the case $v=111, k=11, \lambda=1$), these are the only known parametric restrictions on symmetric 2-designs. The object of this paper is to find  analogous answers to the following question : given a strongly regular graph $G$ and a positive integer $\mu$, when does a quasi-symmetric 2-design with block graph $G$ and defect $\mu$ exist?

The proof of Theorem 1.2 due to Shrikhande uses  p-adic Hilbert symbols,  p-adic invariants for rational equivalence of non-singular symmetric matrices over rational numbers, and a little portion of the Hasse-Minkowski local-global theorem. Following the beautiful little book \cite{Serre} of Serre, we give a brief exposition of this theory in the next section (Section 3). (Another exposition of this theory is  available in Appendix A of \cite{Raghavarao}.) This section continues with a few computational lemmas. Of these, Lemma 3.8 gives a formula for the p-adic invariant of a non-singular symmetric matrix $A$ over $\mathbb{Q}$ in terms of the corresponding invariants of a non-singular principal sub-matrix $B$ and of the \textbf{Schur complement} $A/B$. This lemma, though easy to prove, is likely to be of independent interest. Another important result here is Corollary 3.10 which unearths an unexpected property of rational equivalence. This result may be of wider importance in design theory.

 Section 4 contains the main result and its applications. Let us define an \textbf{integral strongly regular graph}  to be an s.r.g. whose eigenvalues are all integers. Thus, by Theorem 2.3, the block graph of any quasi-symmetric 2-design is an integral s.r.g. We introduce, for any connected integral s.r.g. $G$, the discriminant $\delta (G)$ (taking values in the quotient group $ \mathbb{Q}^{\times}/ \mathbb{Q}^{\Box} $, where $\mathbb{Q}^{\times}$ is the multiplicative group of non-zero rational numbers, and $\mathbb{Q}^{\Box}$ is its subgroup consisting of the squares of non-zero rational numbers) and the p-adic invariant $\epsilon_p(G)$ (taking values in the multiplicative group $\{+1,-1\}$), where $p$ runs over the prime numbers. This is Definition 4.1. In general, these invariants depend on the structure of the s.r.g. $G$, and not merely on its usual parameters. A series of examples illustrating this point appears towards the end of Section 4.

We continue Section 4 by stating and proving the main result (Theorem 4.2) of this paper. It gives parametric restrictions on quasi-symmetric 2-designs solely in terms of the defect of the design, the parameters of a potential block graph $G$ and the invariants $\delta(G)$ and $\epsilon_p(G)$ (p prime).Theorem 4.2 (a) is an analogue of Schutzenberger's Theorem (Theorem 1.1) for quasi-symmetric 2-designs. No such analogue appears to have been known in such generality. Theorem 4.2(b) is an analogue of Shrikhande-Chowla-Ryser Theorem (Theorem 1.2) for quasi-symmetric 2-designs. Similar looking p-adic restrictions (depending only on the parameters of the design, inclusive of the intersection numbers) are available in the literature, see \cite{Bagchi},\cite{Blokhuis}, \cite{Calder1}, \cite{Calder2}, \cite{Calder3} and \cite{Tonchev}.
In general, these results can not be compared with Theorem 4.2(b) since the latter involves the new graph invariants. We explicitly compute these invariants for the complete multipartite graphs, the co-triangular graphs, the non-orthogonality graphs of symplectic spaces over the field of order two and the Steiner graphs (block graphs of 2-designs with $\lambda=1$) hence obtaining new analogues (Corollaries 4.3, 4.5, 4.6 and 4.7) of Theorems 1.1 and 1.2 for quasi-symmetric 2-designs whose block graphs belong to one of these four classes. We  posit a conjecture on the possible parameters of complete multipartite graphs which may occur as the block graphs of q.s. $2$-designs.
We also include a table of small feasible parameters of quasi-symmetric 2-designs with Steiner block graphs.

 As explained in the remarks preceding and following Theorem 4.2, the proof of this theorem is inspired by Shrikhande's paper \cite{SSS2} (published in the year in which this author was busy being born!), and brings its method to a logical conclusion.  Unsurprisingly, Corollary 4.3 of this theorem generalises the main result of \cite{SSS2}, which gives parametric restrictions on affine resolvable 2-designs (i.e., quasi-symmetric 2-designs with complete multi-partite block graphs and $\lambda_1 = 0$).

 \textbf{Acknowledgement.}  We thank Peter Cameron for several fruitful discussions over e-mail pertaining to his construction of quasi-symmetric 2-designs with symplectic block graphs (briefly discussed in Section 4) and regarding the proof of Corollary 2.9 ( which shows
  that the triangular graphs $T_m$ are never the block graphs of quasi-symmetric 2-designs of defect $\mu >1$). We thank B. Sury for acquainting us with Baker's Theorem in \cite{Baker}. A big thank you to Aranya Bagchi, son of the author, for his help with the latex formatting of this paper.

 \section{\textbf{Elementary Restrictions}}

Let $G$ be a finite graph, say on $b$ vertices. The \textbf{adjacency matrix} of $G$ is the $b \times b$ matrix $A$, with its rows and columns indexed by the vertices of $G$, such that, for vertices $x,y$, the $(x,y)$th entry $A(x,y)$ of $A$ is $=1$ if $x,y$ are adjacent in $G$, and $A(x,y)=0$ otherwise. The spectrum $\text{spec}(A)$ (i.e., the multi-set of eigenvalues of $A$, counting multiplicity) is also called the spectrum of $G$, and is denoted $\text{spec}(G)$.

If $G$ is a connected strongly regular graph then $G$ has exactly three distinct eigenvalues,
which we shall denote by $a > \rho > \sigma$, with corresponding multiplicities $1, f, g$. Here $a$ is the degree of $G$ and $f+g+1=b$, the number of vertices of $G$. Since $\text{trace}(A)=0$, we have $a = -f \rho - g \sigma $. We shall refer to $\rho, \sigma, f,g$ as the \textbf{spectral parameters} of the strongly regular graph. One may write the usual parameters of a strongly regular graph in terms of its spectral parameters, and vice versa. Except for the \textbf{conference graphs} (which are the s.r.g.'s with $f=g= (b-1)/2$, eigenvalues $\frac{1}{2}(-1 \pm \sqrt{b})$) on a non-square number of vertices, the eigenvalues of other s.r.g.'s are integers. We shall refer to the  s.r.g.'s with integral eigenvalues as the \textbf{integral strongly regular graphs}.

If $D$ is an incidence system with $v$ points and $b$ blocks, then the \textbf{incidence matrix} $N$ of $D$ is the $v \times b$ matrix, with rows indexed by the points and columns by the blocks of $D$, such that, for a point $x$ and block $B$, the $(x,B)$th entry $N(x,B)$  of $N$ is $=1$ if $ x\in B$, and $N(x,B)=0$ otherwise. For positive integers $n$, $I_n$ and $J_n$ will denote the $n \times n$ identity matrix and the $n \times n$ all-one matrix, respectively.

 If $N$ is the incidence matrix of a 2-design $D$, then $NN^{\prime} = \nu I_v + \lambda J_v$, where $\nu, \lambda$ are the order and balance of $D$. Since $\nu + \lambda v = rk$, it follows that $\text{spec}(NN^{\prime}) = (rk)^1 (\nu)^{v-1}$.

If $D$ is a quasi-symmetric 2-design with usual parameters $b,v,r,k, \lambda$ and intersection numbers $\lambda_1 < \lambda_2$, then the complementary 2-design $\overline{D}$ is also quasi-symmetric, with corresponding parameters given by

\begin{eqnarray} \label{complement} k + \overline{k} = v=\overline{v},\, r + \overline{r} = b = \overline{b},\, \overline{r}-\overline{ \lambda}= r-\lambda,\, \overline{k}-\overline{\lambda}_i = k-\lambda_i \;(i=1,2).\end{eqnarray}
It follows that a complementary pair of quasi-symmetric 2-designs has a common order ($r-\lambda$) and common defect ($\lambda_2-\lambda_1$). Up to isomorphism, they also have the same block graph. The following theorem shows that the parameters of a quasi-symmetric 2-design are determined, up to complementation, by the defect $\mu$ and the spectral parameters of the block graph $G$.
\newline

\textbf{Theorem 2.1}

Let $D$ be a quasi-symmetric 2-design with a connected  block graph $G$ and defect $\mu$. Then the parameters (including intersection numbers) of $D$ are given in terms of the defect $\mu$ and the spectral parameters of $G$ by the following formulae :

\begin{flushleft}
$(\text{a})\; v=f+1,\; b=f+g+1,\; r-\lambda = (\rho - \sigma)\mu,\; k-\lambda_1 = -\sigma \mu,$

$ (\text{b})\; \lambda_2 = \lambda_1 + \mu,\; \lambda_1 =\frac{1}{f+g+1}((f+1)(\lambda + \rho \mu) +g \sigma \mu),$
where $\lambda$ is one of the two roots of the quadratic equation
\begin{equation} \label{quad} x^2 - (f+g+1- 2(\rho - \sigma)\mu)x +(\rho-\sigma)\mu \left((\rho - \sigma)\mu- \frac{f+g+1}{f+1}\right) =0. \end{equation}

(c) Using the second root of this quadratic equation as the value of $\lambda$ in the first two parts of this theorem
  yields the parameters of the complementary design $\overline{D}$.

\end{flushleft}

\textbf{Proof} : Let $A$ and $N$ be the adjacency matrix of $G$, and the incidence matrix of $D$, respectively. We see that $N^{\prime}N = kI_b  + \lambda_2 A + \lambda_1(J_b - I_b - A) = \mu A + (k-\lambda_1)I_b +\lambda_1 J_b$.  Since $\text{spec}(A)=(-f \rho - g \sigma)^1(\rho)^f(\sigma)^g$ and $A$ commutes with $J_b$, it follows that $\text{spec}(N^{\prime}N)= (k + \lambda_1 (b-1) - (f \rho + g \sigma) \mu)^1(k-\lambda_1 + \rho \mu)^f(k-\lambda_1 + \sigma \mu)^g$. But $\text{spec}(NN^{\prime}) = (rk)^1(r-\lambda)^{v-1}$ and $b > v$ (Fisher's inequality) imply that $\text{spec}(N^{\prime}N) = (rk)^1(r-\lambda)^{v-1}(0)^{b-v}$. Comparing these two formulae for the spectrum of $N^{\prime}N$, we get Part (a).
Clearly $\lambda_2 = \lambda_1 + \mu$. Substituting $f+g+1, f+1, \lambda + (\rho - \sigma)\mu, \lambda_1 -\sigma \mu $
for  $b,v,r,k$ (respectively), in the formula $bk=rv$ we get the formula for $\lambda_1$ in Part (b).
Using Equation (\ref{feasible}) and Equation (\ref{complement}), one sees that $\lambda + \overline{\lambda}= b - 2(r-\lambda)= f+g+1-2(\rho -\sigma)\mu$ (which is the negative of the coefficient of $x$ in Equation (\ref{quad})) and $\lambda \overline{\lambda}= (r-\lambda)(r-\lambda-b/v)$, which is the constant term of Equation (\ref{quad}). Thus, $\lambda$ and $\overline{\lambda}$ are the two roots of Equation (\ref{quad}). This completes the proof of Part (b). Applying Parts (a) and (b) of this theorem to $\overline{D}$ in place of $D$, we see that using  $\overline{\lambda}$ in place of $\lambda$ in Parts (a), (b) yields the parameters of $\overline{D}$. This proves Part (c). $\Box$
\newline

We note that Theorem 2.1 places strong restrictions on the possible values of the defect $\mu$ of quasi-symmetric 2-designs with a given block graph. Let us introduce:

\textbf{Definition 2.2} Let $G$ be a connected strongly regular graph, and $\mu$ be a positive integer. We shall say that the pair $(G, \mu)$ is \textbf{feasible} (for a quasi-symmetric design with block graph $G$ and defect $\mu$) if the (complementary pair of)  design parameters (inclusive of the intersection numbers) given by Theorem 2.1 are non-negative integers satisfying the feasibility condition Equation (\ref{feasible}). (In that case, it is easy to see that, with the possible exception of $\lambda_1, \overline{\lambda_1}$,  all these parameters are positive.)
\newline

\textbf{Theorem 2.3}

Let $G$ be a connected strongly regular  graph  with spectral parameters $\rho,\sigma, f, g$, and let $\mu$ be a positive integer. The pair $(G, \mu)$ is feasible for a quasi-symmetric 2-design if and only if the following conditions hold :

\begin{flushleft}

   (a)\; $ \rho \; \mbox{\text{and}}\; \sigma \; \mbox{\text{are integers (i.e.,}}\; G \; \mbox{\text{is an integral s.r.g.),}}$

    (b)\; $- \frac{f+1}{f+g+1}\cdot \frac{f \rho +(g+1)\sigma}{\sigma^2} \leq \mu \leq -\frac{f+1}{2\sigma} $,

    (c)\; $f+1\; \mbox{\text{divides}}\; g(\rho - \sigma) \mu,$

    (d)\; $f+g+1\; \mbox{\text{divides}}\; fg(\rho-\sigma)\mu,\; \mbox{\text{and}}$

    (e)\;$ \mbox{\text{the quantity}}\; \Delta := (f+g+1)(f+g+1 - \frac{4f(\rho - \sigma)\mu}{f+1})\;\mbox{\text{is a perfect square}}. $

\end{flushleft}

(Note that (c) implies that $\Delta$ is an integer.)
\newline

\textbf {Proof} :  Since the characteristic polynomial of any adjacency matrix is a monic integral polynomial, the eigenvalues of any graph are algebraic integers. Thus, (a) holds if and only if $\rho$ and $\sigma$ are rational numbers. But Theorem 2.1(a) shows that rationality of  $\rho , \sigma$ is necessary for the feasibility of $(G,\mu)$. Note that $\Delta$ is the discriminant of the polynomial in Equation (\ref{quad}). Thus, for $\lambda, \overline{\lambda}$ to be rational, it is necessary and sufficient that $\Delta $ is a rational squared. So, once (c) is established, it will follow that (e) also holds. In the rest of the proof, we may assume that $\rho, \sigma$ are integers. Now, the proof of Theorem 2.1 shows that the parameters given there satisfy Equation (\ref{feasible}). Also, the formulae given there show that, (assuming (a)) for all the parameters to be non-negative integers, it is enough to have that $\lambda, \overline{\lambda}$ are integers and $\lambda_1, \overline{\lambda_1}$ are non-negative integers, where over-line denotes (as before) the corresponding complementary parameters.

So, to complete the proof, it suffices to show that  (c) is the necessary and sufficient condition for $\lambda, \overline{\lambda}$ to be integers, (d) is the necessary and sufficient condition for $\lambda_1,\overline{\lambda_1} $ to be integers, and (b) is the necessary and sufficient condition for $\lambda_1, \overline{\lambda_1}$ to be non-negative. Note that, by Theorem 2.1, $\lambda = \gamma \lambda_1 + \delta$, $\overline{\lambda} = \gamma \overline{\lambda_1} + \delta$, and $\lambda, \overline{\lambda}$ are the two roots of the quadratic polynomial (i) $X^2 - \alpha X + \beta$, where,

$$ \alpha =f+g+1- 2(\rho - \sigma)\mu, \; \beta =(\rho-\sigma)\mu((\rho - \sigma)\mu- \frac{f+g+1}{f+1}),$$
$$ \gamma = \frac{f+g+1}{f+1},\; \delta =-\frac{\mu}{f+1}((f+1)\rho +g \sigma).  $$

(Note that $(f+1)\rho +g \sigma = \rho-a < 0$, where $a$ is the largest eigenvalue (degree) of $G$. So we have $\gamma >0, \delta > 0$. This is why the non-negativity of $\lambda_1,\overline{\lambda_1}$ implies positivity of $\lambda, \overline{\lambda}$.) Therefore, substituting $X = \gamma Y + \delta$ in the polynomial (i), we see that  $\lambda_1, \overline{\lambda_1}$ are the two roots of the quadratic polynomial (ii) $Y^2 - \alpha_1Y +\beta_1$, where

$$ \alpha_1 = f+1 + 2 \sigma \mu, \; \beta_1 = \sigma^2 \mu^2 + \frac{f+1}{f+g+1} (f \rho + (g+1)\sigma)\mu. $$

Now observe that, given a monic polynomial of degree two with rational roots, the roots are both integers if and only if all the coefficients are integers, and both the roots are non-negative if and only if the coefficient of its degree one term is $\leq 0$ and the constant term is $\geq 0$. Applying this observation to the polynomials (i) and (ii) completes the proof. $\Box$
\newline

\textbf{Example 0 : Conference graphs.} A conference graph is a strongly regular graph  with spectral parameters $f=g=\frac{q-1}{2}$, $\rho = \frac{1}{2}(-1 + \sqrt{q}), \sigma = \frac{1}{2}(-1 - \sqrt{q})$, where $q \equiv 1\pmod{4}$ is the number of vertices. Applying Theorem 2.3 to such a graph, part (a) shows $q$ must be a perfect square, while part (c) shows that $\frac{q+1}{2}$ divides $\mu$, so that $\mu \geq  \frac{q+1}{2}$. But part (b) shows that we must have $\mu \leq \frac{1}{2}(\sqrt{q}-1)$, contradiction. Thus, conference graphs can never occur as block graphs of quasi-symmetric 2-designs. This was originally observed in \cite{Pawale3}
\newline

\textbf{Example 1 : Complete multi-partite graphs.} The complete multipartite graph $K_{m\times n}$ ($m \geq 2, n\geq 2$)  has $mn$ vertices partitioned into $m$ parts of size $n$ each, where two vertices are adjacent if and only if they are in different parts. In other words,  $K_{m \times n}$ is the complement of $mK_n$ (the disjoint union of $m$ copies of the $n$-vertex complete graph $K_n$). The quasi-symmetric 2-designs with complete multi-partite block graphs are known as \textbf{the strongly resolvable 2-designs}. Recall that, for $n \geq 2$, the $2-(n^2, n, 1)$ designs are known as the affine planes of order $n$. These are strongly resolvable designs of defect $1$ with block graph $K_{n+1 \times n}$.
\newline

\textbf{Corollary 2.4.}

The feasible parameters of strongly resolvable $2$-designs are in bijection with the ordered quadruples $(\alpha, l,l^*,t)$ of non-negative integers such that $\alpha > 0, \,ll^*=\alpha (\alpha -1)$ and $\alpha$ divides $(l+l^*)^2t$. The feasible parameters corresponding to the quadruple $(\alpha, l, l^*,t)$ are given by
\begin{flushleft}
$ n= l+l^*+2\alpha, \; m = \frac{t}{\alpha}n^2+n+1, \; \mu = (n-1)t + \alpha,$
$b=mn,v=n^2((n-1)\frac{t}{\alpha}+1), r= m(\alpha + \ell), \, k = n((n-1)\frac{t}{\alpha}+1)(\alpha + \ell),$
$ \lambda = (\frac{t}{\alpha}n+1)(\alpha + \ell)^2 + \ell, \lambda_1 = n((n-1)\frac{t}{\alpha}+1) \ell, \, \lambda_2 =((n-1)\frac{t}{\alpha}+1) (\alpha + \ell)^2.$
\end{flushleft}
(It easily follows from this result that the only feasible parameters of  strongly resolvable $2$-designs of defect $1$ are those of the affine planes and their complements; further, for each $\mu >1$, there are only finitely many feasible parameters of defect $\mu$.)\newline

\textbf{Proof} : Let the block graph be $G = K_{m \times n}$ ($m\geq 2, n \geq 2$). Its spectral parameters are
$$f=m(n-1), g = m-1, \rho = 0, \sigma = -n .$$
Therefore, in this case, the feasibility conditions of Theorem 2.3 reduce to :
\begin{flushleft}
$(b)\;  mn-m+1 \leq n^2 \mu ,$
\newline
$(c) \; mn-m +1 \; \mbox{\text{divides}} \; (m-n-1)\mu, \, \mbox{\text{and}}$
\newline
$(e) \; \Delta_0 :=(mn-m+1)(mn-m+1 - 4(n-1)\mu) \, \mbox{\text{ is a perfect square.}}$
\end{flushleft}

(Parts (a) and (d) of Theorem 2.3 are automatic. To see that (c) of Theorem 2.3 reduces to (c) above, note that $(m-1)n  \equiv m-n-1 \pmod{mn-m+1}$. The lower bound in part (b) of this theorem follows from part (e) in this case. Indeed, part (e) above implies that $4(n-1)\mu \leq (n-1)m +1$ and hence $m \geq 4\mu-1$, with equality only for $n=2$.)

If $m=4\mu-1$, then $n=2$. Clearly, in this case, the parameters are as in the statement, with $\alpha = 1, l=l^*=0, t =\mu-1$. (These are the parameters of Hadamard $3$-designs.) Therefore, in what follows, we may assume that $m \geq 4\mu$.

Let $(n-1)t$ be the multiple of $n-1$ which is nearest to $\mu$. In case $\mu$ is equidistant between two multiples of $n-1$, we take $(n-1)t$ to be the smaller of them (in which case we have $\mu > (n-1)t$, of course.)

 Suppose, if possible, that $\mu < (n-1)t$. Then $t \geq 1$ and $\mu = (n-1)t-\alpha$, where $1 \leq \alpha < \frac{n-1}{2}$. Hence $n \geq 4$. Since $(n-1)m +1$ divides $(m-n-1)\mu = (m-n-1)((n-1)t-\alpha) = ((n-1)m +1)t -(\alpha m +n^2t -(n+1)\alpha)$, it follows that $(n-1)m+1$ divides $\alpha m +n^2t -(n+1)\alpha$. But, as $t \geq 1$ and $\alpha < \frac{n-1}{2}$, we have $\alpha m +n^2t -(n+1)\alpha > 0$. Therefore $(n-1)m +1 \leq \alpha m +n^2t -(n+1)\alpha$. Hence $\frac{n-1}{2}m < (n-\alpha -1)m < n^2t-(n+1)\alpha < 2(n-1)((n-1)t-\alpha)= 2(n-1)\mu$. (Here, the last inequality holds since $t \geq 1, \alpha < \frac{n-1}{2}$ and $n \geq 4$.) Thus $m < 4 \mu$, a contradiction. Therefore $\mu \geq (n-1)t$, and hence $\mu = (n-1)t + \alpha$, where $0 \leq \alpha \leq \frac{n-1}{2}$.

Since $(n-1)m +1$ divides $(m-n-1)\mu=(m-n-1)((n-1)t+\alpha) = ((n-1)m+1)t +(\alpha m - n^2t -(n+1)\alpha)$, it follows that

\begin{equation} \label{divide} (n-1)m+1 \;\; \text{divides} \;\; \alpha m - n^2t -(n+1)\alpha. \end{equation}

If $\alpha m  < n^2t +(n+1)\alpha$, then (\ref{divide}) implies that $ (n-1)m +1 \leq n^2 t +(n+1)\alpha - \alpha m$, i.e.,
$(n+\alpha-1)m < n^2t + (n+1)\alpha \leq 4(n+\alpha-1)((n-1)t+\alpha)=4(n+\alpha-1)\mu.$ Hence, $m < 4\mu$, a contradiction. Thus, $\alpha m\geq n^2t +(n+1)\alpha$.

If $\alpha m > n^2t +(n+1)\alpha$, then (\ref{divide}) implies that $(n-1)m+1 \leq \alpha m - (n^2t +(n+1)\alpha)\leq \alpha m$, which is absurd since $\alpha <n-1$. Therefore, $\alpha m = n^2t +(n+1)\alpha$. If $\alpha=0$ then this implies $\alpha=t=0$ and hence $\mu=0$, contradiction.
So we have
\begin{equation}\label{formula1} \alpha >0, t \geq 0,\; \mu =(n-1)t+\alpha,\; m= \frac{n^2t}{\alpha}+n+1. \end{equation}

Hence we get $\Delta_0 =(\frac{n}{\alpha})^2((n-1)t+\alpha)^2(n^2-4(n-1)\alpha)$. Since $\Delta_0$ is a square by (e), it follows that $n^2 -4(n-1)\alpha = y^2$ for some integer $y$. Thus, $\alpha (\alpha-1)=\left(\frac{n-2\alpha+y}{2}\right)\left(\frac{n-2 \alpha-y}{2}\right)$. Since $n-2\alpha \pm y$ are integers of the same parity and their product is even, it follows that $l := \frac{n-2\alpha+y}{2}$ and $l^* := \frac{n-2\alpha-y}{2}$ are integers. They satisfy $ll^*=\alpha (\alpha-1)\geq 0, \, l+l^*=n-2\alpha \geq 0$, so that $l,l^*$ are non-negative integers. Since $m=\frac{n^2t}{\alpha}+n+1$, we have $n^2t \equiv 0 \pmod{\alpha}$. Since $n \equiv l+l^* \pmod{\alpha}$, it follows that $(l+l^*)^2 t \equiv 0 \pmod {\alpha}$.
Thus we have
\begin{equation}\label{formula2} l,l^*\geq 0, \; (l+l^*)^2 t \equiv 0 \pmod {\alpha}, \;  ll^*= \alpha (\alpha-1), \; n = l+l^* +2\alpha. \end{equation}

The formulae (\ref{formula1}) and (\ref{formula2}) show that $n,m, \mu$ are given in terms of the quadruple $(\alpha, l,l^*,t)$ as in the statement. Conversely, if $n,m, \mu$ are thus given, with $\alpha >0, \alpha (\alpha-1)=ll^*, (l+l^*)^2t \equiv 0 \pmod{\alpha}$, then it is easy to see that $m,n,\mu$ satisfy (b),(c) and (e), so that the pair $(K_{m\times n}, \mu)$ is feasible. Since $\Delta_0 = (n((n-1)\frac{t}{\alpha} +1)(l-l^*))^2$ by the above computation, it follows from Theorem 2.1 that the design parameters are as  given. Since, in turn, $m,n, \mu$ determine $\alpha, t$ and $l+l^*, ll^*$ (and hence also $l,l^*$ up to a transposition) by the formulae (\ref{formula1}) and (\ref{formula2}), it follows that these formulae give a bijection between the ordered quadruples $(\alpha,l,l^*,t)$ as above and the feasible parameters of strongly resolvable $2$-designs (and, retaining the values of $\alpha,t$ while interchanging $l,l^*$ yields the complementary parameters). $\Box$

\textbf{Example 2: Co-triangular graphs.} For any graph $G$, let $l(G)$ denote the \textbf{line graph} of $G$. Thus, the vertices of $l(G)$ are the edges of $G$, two edges of $G$ are adjacent in $l(G)$ if and only if they meet in one vertex.  The line graphs  $T_n := l(K_n)$ ($n \geq 5$) of the complete graphs $K_n$   are known as the triangular graphs. The co-triangular graph $T_n^*$ is the complement of  $T_n$. Thus, $T_n^*$ may be described as the graph whose vertices are the $\binom{n}{2}$ edges of $K_n$; two edges of $K_n$ are adjacent in $T_n^*$ if and only if they are disjoint. We define a \textbf{Co-triangular 2-design} to be any quasi-symmetric 2-design with a co-triangular block graph.

 If $B$ is a block of a symmetric 2-design $D$, then the 2-design $D_B$ whose blocks are the sets $C\setminus B$, where $C$ runs over the blocks $C \neq B$ of $D$, is known as the \textbf{residual} of $D$ at the block $B$. If $D$ is a \textbf{biplane} (i.e., a symmetric 2-design with $\lambda =2$) with block size $n$, then its residual (at any fixed block) is a $2-(\binom{n-1}{2}, n-2,2)$ co-triangular 2-design with intersection numbers $\lambda_1 =1, \lambda_2 =2$. In \cite{Connor}, Hall and Connor proved that, conversely, any $2-(\binom{n-1}{2}, n-2,2)$ design is a residual of a uniquely determined biplane.
\newline

\textbf{Corollary 2.5} :

(1) The only co-triangular 2-designs of defect $\mu =1$ are the residuals of biplanes and their complements,

(2) For each fixed integer $\mu \geq 2$, there are only finitely many co-triangular 2-designs of defect $\mu$.

(3) The feasible pairs $(T_n^*, \mu)$ are in bijection with the pairs $(\ell, \ell^*)$ of non-negative integers. The bijection is given by the formula $\ell \ell^* = 4 \mu (\mu -1), n = 4 \mu +1 + \ell + \ell^*$.
\newline

\textbf{Proof} : Let the block graph be $T_n^*$. The spectral parameters of $T_n^*$ are
$$f= \frac{n(n-3)}{2}, g= n-1, \rho = 1, \sigma = -(n-3). $$
Therefore, Theorem 2.1 implies that, up to complementation, the co-Steiner designs of defect $1$ are    $2-(\binom{n-1}{2}, n-2,2)$ designs. By the Hall-Connor theorem, these are just the residuals of biplanes. This proves part (1). Part (2) is immediate from part (3) since, given $\mu \geq 2$, the fixed positive integer $4 \mu (\mu-1)$ has only finitely many factors. So it suffices to prove part (3).

 Theorem 2.3 shows that the only feasibility condition is the existence of an integer $x$ such that $(n-4 \mu -1)^2 - 16 \mu (\mu -1) = x^2$. (This is Condition (e) of Theorem 2.3 in this case. This condition implies the upper bound on $\mu$ in part (b) of that theorem. The remaining parts are trivial here.) Hence, $(\frac{n-4 \mu -1 + x}{2})(\frac{n-4 \mu -1 -x}{2}) = 4 \mu (\mu-1)$. As in the proof of Corollary 2.4, the numbers $ \frac{n-4 \mu -1 \pm x}{2} $ are integers. Letting $\ell, \ell^*$ denote these two numbers, we get the result. $\Box$
\newline

\textbf{Example 3 : Symplectic graphs.} Let $d \geq 2$ and let $q$ be a prime power. Take a $(2d)$-dimensional vector space $V$ over the field of order $q$ equipped with a non-degenerate symplectic bilinear form $< \cdot, \cdot > $ (such a form is unique up to linear isomorphisms).
Let $P(V) = PG(2d-1, q)$ be the corresponding projective space. For non-zero vectors $x \in V$, let $[x]$ denote the point in $P(V)$ with homogeneous co-ordinates $x$. The symplectic graph $Sp(2d,q)$ has the points of $PG(2d-1,q)$ as its vertices. Two points $[x], [y]$ are adjacent in $Sp(2d,q)$ if $<x,y> \neq 0$. In short, $Sp(2d, q)$ is the non-orthogonality graph of a symplectic space.
\newpage

\textbf{Corollary 2.6}

Let $q>2$ be a prime power and let $d \geq 2$ be an  integer. Then, a quasi-symmetric $2$-design of defect $\mu$ with block graph $Sp(2d, q)$ is parametrically feasible if and only if  $q(q^{d-1}-1) \equiv 6 \pmod{8}$, $ \mu = (q^d -q +2)/8$, and the pair $(q,d)$ satisfies
$$ \left(\frac{q^d-1}{q-1}\right)^2- q^d\left(\frac{q^{d-1}-1}{q-1}\right)=x^2 $$
for some integer $x$.
\newline

\textbf{Proof} : The spectral parameters of $ Sp(2d,q)$ are
$$\rho = q^{d-1}, \, \sigma = -q^{d-1}, \, f = \frac{\frac{q}{2} (q^{d-1}-1)(q^d +1)}{q-1}, \, g= \frac{\frac{q}{2} (q^{d-1} +1)(q^d -1)}{q-1}. $$

Thus, $f+1 = \frac{1}{2} (q^d-q+2)\frac{q^d-1}{q-1}, \, f+g+1 = \frac{q^{2d}-1}{q-1}$. Therefore, in this case, Condition (c) of Theorem 2.3 simplifies to : $q^d-q+2 | 4q^{d-1}(q-1)$. Since $q >2$, the greatest common divisor between $q^d-q+2$ and $q^{d-1}$ (respectively, between $q^d-q+2$ and $q-1$) is  $1$ or $2$ (respectively $2$ or $1$) according as $q$ is odd or even. Therefore $q^d-q+2$ divides $8 \mu$, say $8 \mu = (q^d-q+2)t$ ($t \geq 1$). Then Condition (e) of Theorem 2.3 simplifies to :

$$ \left(\frac{q^d-1}{q-1}\right)^2- tq^d\left(\frac{q^{d-1}-1}{q-1}\right)=x^2 $$
for some integer $x$. Therefore,
 $$ t \leq \frac{(q^d -1)^2}{q^d(q-1)(q^{d-1}-1)} <2. $$

 Thus $t=1$,  $ \mu = (q^d -q +2)/8$ (so that  $q(q^{d-1}-1) \equiv 6 \pmod{8}$), and the last but one display amounts to the diophantine equation of the statement. Since the remaining conditions of Theorem 2.3 are automatic (with this value of $\mu$), the proof  is complete. $\Box$
\newline

\textbf{Corollary 2.7} : For each fixed integer $d \geq 2$, there are at most finitely many prime powers $q$ for which $Sp(2d,q)$ is the block graph of a quasi-symmetric 2-design.
\newline

\textbf{Proof} : First suppose $d=2$. Then Corollary 2.6 yields $2q+1 = x^2$ for some integer $x$. So, $x$ is odd, say $x=2y+1$. Then $q=2y(y+1)$. Thus $q$ is (even, and therefore) a power of $2$. Hence both $y$ and $y+1$ are powers of $2$, so that $y=1$. Hence $q=4$. But the pair $(d,q)=(2,4)$ fails the congruence condition of Corollary 2.6. Thus, there is no prime power $q$ for which $SP(4,q)$ is the block graph of a q.s. $2$-design.

Next, let $d \geq 3$. A theorem of Alan Baker from \cite{Baker} says that if $f$ is a single variable polynomial over integers which has at least three distinct simple roots, then there can only be finitely many integers $q$ such that $f(q)$ is a perfect square. Consider the polynomial $f_d(X) := (\frac{X^d-1}{X-1})^2 - X^d(\frac{X^{d-1}-1}{X-1})$. It may be shown that, for $d \geq 3$, $f_d$ satisfies Baker's hypothesis. Therefore, for $d \geq 3$, this result follows from Baker's theorem and Corollary 2.6. $\Box$

It also seems likely that for each fixed prime power $q \geq 3$, the graph $Sp(2d,q)$ is a block graph for at most finitely many values of $d$. Probably this can be deduced from (Corollary 2.6 and) the ABC conjecture.
\newline

\textbf{Example 4 : Steiner graphs.}  Recall that a \textbf{Steiner 2-design} is a 2-design with balance $\lambda = 1$. Note that any Steiner 2-design is automatically quasi-symmetric of defect $1$.  For integers $m > n \geq 2$ such that $n$ divides $m(m-1)$, we define the \textbf{Steiner graph} $ S_n(m)$ to be the block graph of a Steiner 2-design with parameters
$$ b= \frac{m}{n}(mn-m+1), v= mn-m+1, r=m,k=n, \lambda=1.$$
Note that, when the parameters $m,n$ are large, there usually are many non-isomorphic Steiner graphs, all of them designated $S_n(m)$.
We define a \textbf{multi-Steiner 2-design} to be any quasi-symmetric 2-design with a Steiner block graph.
\newline

\textbf{Corollary 2.8}

(1) The only multi-Steiner 2-designs of defect $\mu=1$ are the Steiner 2-designs and their complements,

(2)  For each fixed pair of  integers $\mu \geq 2, n\geq 2$, there are only finitely many multi-Steiner 2-designs of defect $\mu$ with block graph $S_n(\cdot)$.
\newline

\textbf{Proof} : (1) is immediate from the definition of $S_n(m)$ and Theorem 2.1.
By Theorem 2.1, the spectral parameters of $ S_n(m)$ are :
$$ f= m(n-1),\, g= m(m-n +1) -1 - \frac{m(m-1)}{n}, \,  \rho = m-n-1, \, \sigma =-n. $$

Note that, by Theorem 2.3, the only feasibility requirement for the existence of a quasi-symmetric 2-design of defect $\mu$ with block graph $S_n(m)$ is (apart from $n|m(m-1)$) that there is an integer $x$ such that $ (mn-m +1-2n\mu)^2 -4n^2 \mu (\mu -1)=x^2 $.  (In this case, the upper bound in (b) of  Theorem 3.2 follows from this requirement, and the remaining parts are trivial.) We have $(\frac{mn-m+1 -2n \mu +x}{2})(\frac{mn-m+1-2n \mu-x}{2})= n^2 \mu (\mu-1)$. As in the proof of Corollary 2.4, the numbers  $\frac{mn-m+1 -2n \mu \pm x}{2} $ are  integers. Thus, for fixed values of $n \geq 2, \mu \geq 2$,  $\frac{mn-m+1 -2n \mu \pm x}{2} $ are among the finitely many factors of the fixed positive integer $n^2 \mu (\mu -1)$, so that there are only finitely many feasible values of $m$. This proves part (2). $\Box$

Note that the symplectic graph $Sp(2d,2)$ has the same parameters as a Steiner graph $S_n(m)$ with $n=2^{d-1}$, $m = 2^d +1$. Thus, Corollary 2.8 applies to these graphs as well. (More generally, it applies to all \textbf{pseudo-Steiner graphs} : s.r.g.'s having the same parameters as Steiner graphs.) This is why we left out the case $q=2$ in Corollary 2.6.

A $t-(v,k,\lambda)$ design is an incidence system with $v$ points, $k$ points per block, and $\lambda$ blocks containing any $t$ distinct points. Any $t$-design is an $s$-design for each $s$ in the range $0 \leq s \leq t$.  Ray- Chaudhuri and Wilson generalized (\cite{Chaudhuri}) Fisher's inequality to prove that the parameters of any (2s)-design with $v \geq k+s$ satisfy $b \geq \binom{v}{s}.$ The $(2s)$-designs attaining this bound are known as the tight $(2s)$-designs. This generalizes the notion of symmetric $2$-designs (the case $s=1$). Ito in \cite{Ito} and Bremner in \cite{Bremner} proved that, up to complementation, the only tight $4$-design with $4 \leq k \leq v-4$ is the famous $4-(23, 7,1)$ design of Witt. Its block graph is a sporadic strongly regular graph with spectral parameters $f=22, g=230, \rho = 25 , \sigma = -3$. On the other hand, Cameron has proved (\cite{Cameron1}) that any quasi-symmetric $2$-design with a connected block graph satisfy $b \leq \binom{v}{2}$, and equality holds only for $4$-designs. The following Theorem is an immediate consequence of these result.
\newline

\textbf{Corollary 2.9}:

There is no  quasi-symmetric 2-design of defect $\mu \geq 2$ whose block graph is a triangular graph $T_m := S_2(m-1)$ ($m \geq 5$).
\newline

\textbf{Proof} : $T_m$ has spectral parameters
$$ \rho = m-4,\, \sigma = -2,\, f = m-1, \, g= \binom{m-1}{2}-1.$$
 Therefore, by Theorem 2.1, any such design would satisfy $b = \binom{v}{2}$, $4 \leq k \leq v-4$. Therefore, by Cameron's Theorem, the design would be a non-trivial tight $4$-design with block graph $T_m$. But by the classification of Ito and Bremner, there is no such design. $\Box$

   Note that this theorem rules out infinitely many feasible parameters which survive Corollary 4.7 below. Cameron has  pointed out that the upper and lower bounds  on $b$ (for quasi-symmetric $2$-designs and non-trivial $4$-designs, respectively) mentioned above  are very special cases of Theorem 5.21 in Delsarte's Thesis.

\section{\textbf{Hilbert symbols and rational equivalence}}
 $\mathbb{Q}^{\times}$ will denote the multiplicative group of non-zero elements in the field $\mathbb{Q}$ of rational numbers. We denote by $\mathbb{Q}^{\Box}$  the subgroup of $\mathbb{Q}^{\times}$ consisting of the squares of non-zero rationals. For $x,y \in \mathbb{Q}^{\times}$, we write $x \equiv y \pmod {\mathbb{Q}^{\Box}}$ if $xy^{-1} \in \mathbb{Q}^{\Box}$.  Also, for $\alpha \in  \mathbb{Q}^{\times}/ \mathbb{Q}^{\Box}$ and $x \in  \mathbb{Q}^{\times}$ we sometimes write $\alpha = x \pmod {\mathbb{Q}^{\Box}}$ to indicate that $\alpha$ is the image of $x$ under the quotient map  $\mathbb{Q}^{\times} \rightarrow \mathbb{Q}^{\times}/ \mathbb{Q}^{\Box}$.

Through out this section, $p$ is a prime number, fixed but arbitrary. An element $u$ of  $\mathbb{Q}^{\times}$ is said to be a \textbf{p-adic unit} if $p$ does not divide the numerator and denominator of $u$ in its reduced form. The p-adic valuation $v_p(x)$ of $x \in   \mathbb{Q}^{\times}$ is the unique integer $m $ such that $p^{-m}x$ is a p-adic unit. The p-adic norm of $x$ is defined by $\|x\|_p=p^{-v_p(x)}$. This is extended to $\mathbb{Q}$ by setting $\| 0 \|_p = 0$. It is easy to see that $\| \cdot \|_p$ is a field norm on $\mathbb{Q}$. The field $\mathbb{Q}_p$ of p-adic numbers is defined to be the completion of $\mathbb{Q}$ under this norm. Thus, $\mathbb{Q}$ is a subfield of $\mathbb{Q}_p$.

The \textbf{ p-adic Hilbert symbol} is the function from $\mathbb{Q}^{\times} \times \mathbb{Q}^{\times}$ to $\{ +1, -1 \}$  defined as follows. For $a, b \in \mathbb{Q}^{\times}$, $(a,b)_p = +1$ if the equation $ax^2+by^2 = z^2$ has a solution $(x,y,z) \not = (0,0,0)$ in  $\mathbb{Q}_p$, and $(a,b)_p = -1$ otherwise. This symbol has the following important properties (see \cite{Serre}, Chapter III) :-

\begin{flushleft}
(H1)  $\mathbb{Q}^{\Box}$-invariance : For $a,a^{\prime},b,b^{\prime} \in \mathbb{Q}^{\times}$, if $a^{\prime} \equiv a \pmod {\mathbb{Q}^{\Box}}$ and $b^{\prime} \equiv b \pmod {\mathbb{Q}^{\Box}}$ then $(a^{\prime}, b^{\prime})_p = (a,b)_p$.
\newline
(H2) Symmetry : For $a,b \in  \mathbb{Q}^{\times}$, $(b,a)_p=(a,b)_p$.
\newline
(H3) Bilinearity : For $a,b,c \in  \mathbb{Q}^{\times}$, $(ab,c)_p=(a,c)_p(b,c)_p, \, (a,bc)_p= (a,b)_p(a,c)_p$.
\newline
(H4) Special identities :  For $a,b \in  \mathbb{Q}^{\times}$ with $b \neq 1$, $(a,-a)_p= 1 =(-a,a)_p, \, (b,1-b)_p=1=(1-b,b)_p$.
\newline
(H5) Formulae : For p-adic units $u,v \in\mathbb{Q}^{\times}$, we have,
     $$ (u,p)_p = \begin{cases} (-1)^{\omega(u)} \;\text{if}\; p=2, \\ (\frac{u}{p}) \;\text{if}\; p \neq 2.
      \end{cases} $$
      $$ (u,v)_p = \begin{cases} (-1)^{\epsilon (u) \epsilon (v)} \; \text{if}\: p=2,\\ 1 \; \text{if} \; p \neq 2.
      \end{cases} $$
\end{flushleft}

Here, for odd $p$, $(\frac{\cdot}{p})$ is the Legendre symbol : for p-adic units $u$, $(\frac{u}{p})= +1$ if $u$ is a square modulo $p$, and $=-1$ otherwise. For $2$-adic units $u$, $\omega(u)=0$ if $u \equiv \pm 1 \pmod{8}$ and $\omega(u)=1$ if $u \equiv \pm 3 \pmod{8}$; $\epsilon(u)=0$ if $u \equiv +1 \pmod{4}$ and $\epsilon(u)=1$ if $u \equiv -1 \pmod{4}$.

Remark on terminology. The use of the word `bilinearity' to describe (H3) requires an explanation. Note that the target $\{+1, -1\}$ of the Hilbert symbols is a field with ordinary multiplication as field addition; the field multiplication is determined by the requirement that $+1$ is the additive identity and $-1$ the multiplicative identity in this field. Since the quotient group $\mathbb{Q}^{\times}/ \mathbb{Q}^{\Box}$
is a multiplicative elementary abelian $2$-group, it may be viewed as a vector space over this field. The properties (H1)--H(3) say that the p-adic Hilbert symbol descends to a well defined symmetric bilinear form on this vector space. Another important property of this symbol is the non-degeneracy of this bilinear form. That is, an element $x$ of $\mathbb{Q}^{\times}$ satisfies $(x,y)_p=1$ for all $y \in  \mathbb{Q}^{\times}$
(if and) only if  $x \in  \mathbb{Q}^{\Box}$. Non-degeneracy will not be of importance to us.

Since  $\mathbb{Q}^{\times}$ is generated by the p-adic units together with $p$, the value of $(x,y)_p$ may be calculated using the properties (H) for any given elements $x,y$ of  $\mathbb{Q}^{\times}$. Through the rest of this article, we shall use (H) without further mention.

We shall say that a solution of a homogeneous quadratic equation (in several variables) is non-trivial if at least one of the co-ordinates of the solution is non-zero. The following result (Corollary 1 in \cite{Serre}, Chapter IV) is a baby version of the Hasse-Minkowski local-global theorem.
\newline

\textbf{Lemma 3.1}

For $a,b \in \mathbb{Q}^{\times}$, the following two conditions are equivalent :
(1) The equation $ax^2 + by^2=z^2$ has a non trivial solution in rationals (equivalently, in integers),
(2) $(a,b)_p=1$ for all primes $p$.

In other words, the equation in (1) has a non-trivial solution in rationals if and only if it has a non-trivial solution in every
 $\mathbb{Q}_p$. (It is usual to add here the requirement that the equation is solvable in reals as well. But, the Hilbert product formula (Theorem 3 in \cite{Serre}, Chapter III) shows that if $(a,b)_p=1$ for all primes $p$ then this equation is automatically solvable in reals; i.e., in that case, both of  $a,b$ can not be negative.)

Combining standard arguments from elementary number theory with the case $b=-1$ of Lemma 3.1, we get :
\newline

\textbf{Corollary 3.2}

Let $n$ be a positive integer. Then $n$ is a sum of (at most) two perfect squares if and only if $(-1,n)_p=1$ for all primes p.

Now we recall that two $n \times n$ symmetric matrices $A,B$ over $\mathbb{Q}$ are said to be \textbf{rationally equivalent} (in symbols $A \sim B$) if there is a non-singular $n \times n$ matrix $X$ over $\mathbb{Q}$ such that $B=X^\prime AX$. (Here $X^\prime$ is the transpose of
$X$.) Clearly, rational equivalence is an equivalence relation on the space of all $n \times n$ symmetric matrices over $\mathbb{Q}$. Note that, if $A \sim B$ and $A$ is non-singular then so is $B$.

Let $D = \text{diag}(d_1, \cdots, d_n)$ be a non-singular diagonal matrix over  $\mathbb{Q}$. Then \textbf{the p-adic invariant} $\epsilon_p (D)$ of $D$ is defined by
$$  \epsilon_p (D):= \prod_{1 \leq i < j \leq n} (d_i, d_j)_p. $$
The next lemma is essentially Theorem 5 in \cite{Serre}, Chapter IV.

\textbf{Lemma 3.3}

Let $D_1$ and $D_2$ be two  non-singular diagonal matrices over  $\mathbb{Q}$. If $D_1 \sim D_2$ then $ \epsilon_p (D_1) =  \epsilon_p (D_2)$.

Another basic observation is :
\newline

\textbf{Lemma 3.4}

 Every non-singular symmetric matrix over  $\mathbb{Q}$ is rationally equivalent to a (nonsingular) diagonal matrix over  $\mathbb{Q}$.
\newline

\textbf{Proof} : Let $A$ be an $n\times n$ non-singular symmetric matrix over  $\mathbb{Q}$. Define the non-degenerate symmetric bilinear form
$(\cdot, \cdot)$ on $\mathbb{Q}^n$ by $(x,y)= x^\prime Ay$. Using the usual Gram-Schmidt algorithm, any given basis of  $\mathbb{Q}^n$ can be orthogonalised (not orthonormalised : the normalisation is generally impossible over $\mathbb{Q}$) with respect to this bilinear form. This process yields a basis $\{ x_1, \ldots, x_n \}$ such that $x_i^\prime A x_j =0$ for all $i \neq j$ ($1 \leq i,j \leq n$). Let $d_i = x_i^\prime A x_i$, $1 \leq i \leq n$. Let $D= \text{diag} (d_1, \ldots, d_n)$. Let $X$ be the $n \times n$ (non-singular) matrix over $\mathbb{Q}$
whose columns are the vectors $ x_1, \ldots, x_n$. Then $D =X^{\prime} AX $. $ \Box$

Now, for any  $n\times n$ non-singular symmetric matrix $A$ over  $\mathbb{Q}$, \textbf{the p-adic Hasse invariant} of $A$ is defined by
$$\epsilon_p (A) := \epsilon_p(D), $$
where $D$ is any ( non-singular) diagonal matrix over  $\mathbb{Q}$ such that $D \sim A$. Lemma 3.3 and 3.4 show that this is well defined : such a matrix $D$ exists and $\epsilon_p (A)$ is independent of the choice of $D$. Lemma 3.3 also shows that it is indeed an invariant for rational equivalence.The Hasse invariants were introduced by Hasse in 1923, building on previous work of Minkowski. 
\newline

\textbf{ Theorem  3.5}

Let $A_1, A_2$ be two $n\times n$ non-singular symmetric  matrices over   $\mathbb{Q}$. If $A_1$ and $A_2$ are rationally equivalent, then $\det(A_1) \equiv \det(A_2) \pmod {\mathbb{Q}^{\Box}}$ and  $\epsilon_p(A_1)=\epsilon_p(A_2)$.

This completes our mini-survey of rational equivalence and the p-adic invariant. For more on this topic, \cite{Serre} is a perfect source. The next few results in this section may be new. We have failed to locate them in the available literature.
\newline

\textbf{Lemma 3.6}

Let $a,b \in \mathbb{Q}$ be such that $a \neq 0$ and $a+bn \neq 0$. Then,
$$ \epsilon_p(aI_n + bJ_n) = (-1,a)_p^{\binom{n-1}{2}} (a, a+bn)_p^{n-1} (a(a+bn), n)_p. $$

\textbf{Proof} : Note that our hypotheses on $a,b$ are necessary for the matrix $aI_n + bJ_n$ to be non-singular.
Let $X_n$ be the $n \times n$ matrix given by

$$ X_n(i,j)= \begin{cases} 1 \; \text{if} \; j=1, 1 \leq i \leq n, \\
-1 \; \text{if} \; 1 \leq i < j \leq n, \\
j-1 \; \text{if} \; 1 <  i=j \leq n, \\
0 \; \text{if} \; 1 <  j < i \leq n. \end{cases} $$

A calculation yields $X_n^{\prime}X_n = D_1, \; X_n^{\prime} J_n X_n = D_2$, where $D_1 = \text{diag}( n, j(j-1) : 1 < j \leq n)$ and $D_2
= \text{diag}(n^2, 0,0, \ldots, 0)$. Therefore $X_n^{\prime} (aI_n + b J_n) X_n =D$ where $D = a D_1 + b D_2$. Since $D_1$ is clearly non-singular, it follows that $X_n$ is non-singular and $aI_n + bJ_n \sim D$. Hence $\epsilon_p (aI_n + b J_n)= \epsilon_p (D)$ by Theorem 3.5.
We now compute

\begin{flushleft}
$ \epsilon_p(D)= \prod_{1< j \leq n}(n(a+bn), aj(j-1))_p \cdot \prod_{1 < i < j \leq n}(ai(i-1), aj(j-1))_p .$
\newline

But,
$  \prod_{1< j \leq n}(n(a+bn), aj(j-1))_p =  (n(a+bn), a)_p^{n-1} \prod_{1< j \leq n}(n(a+bn), j(j-1))_p. $
\newline

And, as $(a,a)_p = (-1,a)_p(-a, a)_p = (-1,a)_p$,
\newline

$   \prod_{1 < i < j \leq n} (ai(i-1), aj(j-1))_p = (-1,a)_p^{\binom{n-1}{2}}\prod_{1 <j \leq n}(a, j(j-1))_p^{n-2}\prod_{1 < i < j \leq n}( i(i-1), j(j-1))_p . $
\end{flushleft}

Also,

\begin{eqnarray*}  \prod_{1< j \leq n}(n(a+bn), j(j-1))_p & = & \prod_{1< j \leq n}(n(a+bn), j)_p (n(a+bn), j-1)_p \\
                                                           & = & (n(a+bn),n)_p = (-(a+bn), n)_p,\end{eqnarray*}
since this last product is telescoping (remember : the Hilbert symbol is $\pm 1$-valued). Similarly,

$$\prod_{1 <j \leq n}(a, j(j-1))_p = (a, n)_p.  $$

Thus we get

\begin{eqnarray*} \epsilon_p(a I_n + b J_n) & = & (n(a+bn), a)_p^{n-1} (-1,a)_p^{\binom{n-1}{2}}(a,n)_p^n (-(a+bn),n)_p e_n \\
& = & (a+bn,a)_p^{n-1} (-1,a)_p^{\binom{n-1}{2}} (-a(a+bn), n)_p e_n, \end{eqnarray*}

where $e_n =  \prod_{1 < i < j \leq n}( i(i-1), j(j-1))_p $. Therefore, to complete the proof, it suffices to show that $e_n = (-1,n)_p$ for all $n$. Vacuously, $e_1 = 1 =(-1,1)_p$. Also,
$$ e_n e_{n+1} = \prod_{1 < i \leq n}(i(i-1), n(n+1))_p = (n, n(n+1))_p . $$

But, $(n, n(n+1))_p = (-1,n(n+1))_p(-n,n)_p(-n, 1-(-n))_p = (-1, n(n+1))_p$. Thus, $ e_1= (-1,1)_p, \, e_ne_{n+1} = (-1,n)_p(-1, n+1)_p$. Hence, by induction on $n$, we get $e_n = (-1,n)_p$. $\Box$
\newline

\textbf{Lemma 3.7}

Let $A,B$ be non-singular symmetric matrices over $\mathbb{Q}$, not necessarily of the same order. Then,
$$\epsilon_p(A \oplus B) =  \epsilon_p(A)\epsilon_p(B) (\det(A), \det(B))_p. $$

\textbf{Proof} : In view of the definition of $\epsilon_p(\cdot)$, we may assume without loss of generality that both $A$ and $B$ are diagonal matrices, so that $A \oplus B$ is also diagonal. In this case, the result follows from the definition of  $\epsilon_p(\cdot)$ for diagonals and the bilinearity of Hilbert symbols. $\Box$
\newline
Let $A$ be a symmetric matrix over   $\mathbb{Q}$ and let $B$ be a non-singular principal sub-matrix of $A$. Then the \textbf{{\large Schur complement}} $A/B$ of $B$ in $A$ is defined as follows.\newline
 Without loss of generality, we may assume that
$$ A= \left( \begin{array}{cc}
B & C \\ C^{\prime} & D \end{array} \right),
 \mbox{ where $D$ is also symmetric.}$$
 $$ \mbox{  Then, } A/B := D-C^{\prime}B^{-1}C. $$

The first part of the following lemma is one among several little gems due to I. Schur, each of which is known as Schur's Lemma. The second part may be new.
\newline

\textbf{Lemma 3.8}

Let $A$ be a symmetric matrix over  $\mathbb{Q}$ and let $B$ be a non-singular principal sub-matrix of $A$. Then,

(a)  $\det (A/B) = \det (A)/ \det(B) $ and hence $A$ is non-singular \newline if and only if $A/B$ is non-singular.

(b)  If, further, $A/B$ is non-singular, then
$$ \epsilon_p(A) =  \epsilon_p(B)\epsilon_p(A/B) (\det(B), \det(A/B))_p. $$

{\textbf{Proof} : Let $A$ be given by the $2 \times 2$ block matrix in the definition of Schur complement. Let
$$  X := \left( \begin{array}{cc}
I & -B^{-1}C \\ 0 & I \end{array} \right).$$

Since $X$ is a block triangular matrix  with identities as diagonal blocks, we have $\det (X)=1$, and hence $X$ is non-singular. A computation shows that
$$ X^{\prime}AX = B \oplus (A/B). $$

Since $\det (X)=1$, this proves Part (a) and shows that $A \sim B \oplus (A/B) $. Hence by Theorem 3.5, if $A/B$ is also non-singular, then
 $  \epsilon_p(A) =  \epsilon_p(B \oplus(A/B)) $. Therefore Lemma 3.7 completes the proof of Part (b). $\Box$
\newline

 In the next lemma, for any matrix $A$ over  $\mathbb{Q}$, ${\cal C}(A)$ will denote the column space of $A$ over  $\mathbb{Q}$. It is the
  $\mathbb{Q}$-vector space spanned by the columns of  $A$.
\newpage

  \textbf{Lemma 3.9}

  Let $E_1, E_2$ be two $m \times n$ matrices of rank $n$ over  $\mathbb{Q}$. Suppose ${\cal C}(E_1)= {\cal C}(E_2)$. Then
   $E_1^{\prime}E_1$ and $E_2^{\prime}E_2$ are rationally equivalent.
\newline

  \textbf{Proof} : Let  ${\cal C}(E_1)= V ={\cal C}(E_2)$. Let $\{x_1, \ldots, x_n \}$ and  $\{y_1, \ldots, y_n \}$ be the set of columns of $E_1$ and of $E_2$, respectively. These two sets are two bases of the  $\mathbb{Q}$-vector space $V$. Let $A=((a_{ij}))_{1 \leq i,j \leq n}$ be the transition matrix between these two bases. That is, $a_{ij} \in  \mathbb{Q}$ are determined by the equations  $\sum_{i=1}^n a_{ij}x_i = y_j$ for $1 \leq j \leq n$. Thus, $A$ is a non-singular matrix over  $\mathbb{Q}$, and $E_2 = E_1A$. Therefore $ E_2^{\prime}E_2 = A^{\prime}( E_1^{\prime}E_1)A$.
  Hence $ E_2^{\prime}E_2 \sim  E_1^{\prime}E_1$. $\Box$
\newline

\textbf{Corollary 3.10}

Let $E$ be an $m \times n$ matrix of rank $g$ over  $\mathbb{Q}$. Then all the $g \times g$ non-singular principal sub-matrices of $E^{\prime}E$ are rationally equivalent.
\newline

\textbf{Proof} : Let $B_1, B_2$ be two  $g \times g$ non-singular principal sub-matrices of $E^{\prime}E$. So there are rank $g$ sub-matrices $E_1, E_2$ of $E$, both of order $m \times g$, such that $B_i = E_i^{\prime}E_i$ ($i=1,2)$. We have $ {\cal C}(E_1)= {\cal C}(E)= {\cal C}(E_2)$. Hence Lemma 3.9 implies $B_1 \sim B_2$. $\Box$

\section {\textbf{ The main result and applications.}}

If $N$ is the $v\times v$ incidence matrix of a symmetric 2-design on $v$ points with balance $\lambda$ and order $\nu$, then $k^2 \nu^{v-1} = \det(NN^{\prime})=(\det N)^2$ implies $N$ is non-singular and $\nu^{v-1} \equiv 1 \pmod{\mathbb{Q}^{\Box}}$, proving Theorem 1.1. Also, $NN^{\prime} = \nu I_v + \lambda J_v$ implies that $ \nu I_v + \lambda J_v \sim I_v$, and hence, by Theorem 3.5, $\epsilon_p(  \nu I_v + \lambda J_v)=1$. But, since $\nu + \lambda v = k^2 \equiv 1\pmod{\mathbb{Q}^{\Box}}$,  Lemma 3.6 yields that, when $v$ is odd (so that $\binom{v-1}{2} \equiv \frac{v-1}{2} \pmod{2}$),  $\epsilon_p(  \nu I_v + \lambda J_v)=  (\nu, (-1)^{(v-1)/2} v)_p$. But $(\nu, \lambda v)_p = (\nu, k^2-\nu)_p = (\frac{\nu}{k^2}, 1- \frac{\nu}{k^2})_p=1$ and hence $(\nu,v)_p=(\nu,\lambda)_p$. Thus, when $v$ is odd, we must have $ (\nu, (-1)^{(v-1)/2} \lambda)_p = \epsilon_p( \nu I_v + \lambda J_v)= 1 \, \forall p$ for the existence of a symmetric 2-design with these parameters. In view of Lemma 3.1, this proves Theorem 1.2. (It may be instructive to compare this short proof with the proof given in Chapter 12 of \cite{Raghavarao}.)

Clearly these arguments work since the incidence matrix of a symmetric 2-design is a non-singular matrix. The main idea of
\cite{SSS2} in extending such arguments to a class of quasi-symmetric 2-designs was to use a construction of Connor to embed the $v \times b$ incidence matrix as a sub-matrix of a suitable non-singular $b \times b$ matrix $M$ over $\mathbb{Q}$, and then apply the Hasse-Minkowski theory to $MM^{\prime}$.

In Theorem 4.2  below, we show that the technique of \cite{SSS2} works in the generality of all quasi-symmetric 2-designs without any need for Connor's special construction. Indeed, $M$ may be taken to be an arbitrary  non-singular $b \times b$ matrix  over $\mathbb{Q}$ containing the incidence matrix as a sub-matrix (see the remark after Theorem 4.2). To state and prove our result, we need the following definition.
\newline

\textbf{Definition 4.1}: Let $G$ be a connected integral strongly regular graph. Let $E$ be the orthogonal projection onto the negative eigenspace of $G$ (i.e., the eigenspace corresponding to the negative eigenvalue $\sigma$ of $G$). In other words, $E$ is the minimal idempotent of rank $g$ in the Bose-Mesner algebra of $G$. Since $G$ is an integral s.r.g., it follows that $E$ is a matrix over $  \mathbb{Q}$. We define the discriminant $\delta (G) \in \mathbb{Q}^{\times}/ \mathbb{Q}^{\Box}$ and the p-adic invariant $\epsilon_p (G) = \pm 1$ ($p$ any prime number) of $G$ as follows :

$$ \delta (G)= \det(E_0)\pmod{\mathbb{Q}^{\Box}}, \; \epsilon_p (G) = \epsilon_p (E_0) $$
where $E_0$ is any $g \times g$ non-singular principal sub-matrix of $E$. Since $E$ is of rank $g$ and $E^{\prime}E = E$, Corollary 3.10 implies that this definition is independent of the choice of $E_0$.

The next theorem is the main result of this paper. In view of Theorem 2.1, it is natural to present it in terms of the block graph and the defect $\mu$ of the design, or equivalently (since $\nu = (\rho - \sigma )  \mu$ by Theorem 2.1), in terms of the graph parameters and  the order $\nu$ of the design. This result looks neater in terms of the order.
\newline

\textbf{Theorem 4.2}

Let $G$ be a connected integral strongly regular graph with spectral parameters $\rho, \sigma, f,g$, and let $\nu$ be a positive integer. Then, for the existence of a quasi-symmetric 2-design of order $\nu$ with block graph $G$, the following conditions are necessary :

$$(a) \,  \nu^f \equiv (f+1)(f+g+1) \delta (G) \pmod{\mathbb{Q}^{\Box}}, $$
and, for all prime numbers $p$,
$$(b) \, (-1, \nu)_p^{\binom{f}{2}} (\nu, f+1)_p = (f+g+1, -f-1)_p (-(f+1)(f+g+1), \delta(G))_p \epsilon_p(G). $$

\textbf{Proof} : Let $E_0$ be as in Definition 4.1. Without loss of generality, we may assume that $E_0$ is the  $g \times g$ principal sub-matrix of $E$ in its top left corner. That is, the rows (and columns) of $E_0$ correspond to the first $g$ rows (respectively columns) of $E$. Let
$N$ be the $v \times b$ incidence matrix of a quasi-symmetric 2-design of order $\nu$ with block graph $G$. Consider the $b \times b$ matrix $M$ given by

$$M=  \left( \begin{array}{c} N\\ X \end{array} \right), \text{where}\; X= \left( \begin{array}{cc} I_g, & 0_{g \times v} \end{array} \right). $$

(Recall : $g=b-v$ by Theorem 2.1.) Note that $NN^{\prime} = \nu I_v + \lambda J_v$ is non-singular, and we have
$$ MM^{\prime}=  \left( \begin{array}{cc}
NN^{\prime} & NX^{\prime} \\ XN^{\prime} & XX^{\prime} \end{array} \right)$$

Hence the Schur complement of $NN^{\prime}$ in $MM^{\prime}$ is $MM^{\prime}/N N^{\prime} = X(I_b - N^{\prime}(NN^{\prime})^{-1} N)X^{\prime}$. But $ I_b -N^{\prime} (NN^{\prime})^{-1} N$ is the orthogonal projection onto the kernel of $N^{\prime}N$. (This is true of any matrix $N$ such that $NN^{\prime}$ is non-singular.) Also, the proof of Theorem 2.1 shows that the kernel of  $N^{\prime}N$ is precisely the negative eigen-space of $G$. Thus, $ I_b - N^{\prime}(NN^{\prime})^{-1} N = E$, the minimal idempotent of rank $g$ in the Bose-Mesner algebra of $G$. Thus $MM^{\prime}/NN^{\prime} = XEX^{\prime}$. Also, since $E_0$ is the $g \times g$ sub-matrix in the top left corner of $E$, our choice of $X$ implies that $  XEX^{\prime} = E_0$. Therefore, $MM^{\prime}/NN^{\prime} =E_0$.

Since both $E_0$ and $NN^{\prime}$ are non-singular, Lemma 3.8 implies that $MM^{\prime}$ (and hence also $M$) is non-singular, so $\det(MM^{\prime})=(\det M)^2 \equiv 1  \pmod{\mathbb{Q}^{\Box}}$, and $MM^{\prime} \sim I_b$. Thus by Theorem 3.5, $\epsilon_p(MM^{\prime})=1$. Since $r/k = b/v$, so that $rk \equiv bv   \pmod{\mathbb{Q}^{\Box}}$,  we get $\det (NN^{\prime})  = rk \nu^{v-1} \equiv bv \nu ^{v-1} \pmod{\mathbb{Q}^{\Box}}$. Also, $\det(E_0) \equiv \delta (G)  \pmod{\mathbb{Q}^{\Box}}$ by Definition 4.1.
Hence Lemma 3.8(a) implies that $\nu^{v-1} bv \delta (G) \equiv \det (MM^{\prime}) \equiv 1  \pmod{\mathbb{Q}^{\Box}}$. Thus, $\nu^{v-1} \equiv
bv \delta (G)  \pmod{\mathbb{Q}^{\Box}}$. Since $v=f+1, b= f+g+1$ by Theorem 2.1, this proves the first part.
\newline
Since $\nu + \lambda v \equiv bv  \pmod{\mathbb{Q}^{\Box}}$ and $\nu^{v-1} \equiv bv \delta (G)  \pmod{\mathbb{Q}^{\Box}}$, Lemma 3.6 implies

\begin{eqnarray*} \epsilon_p (NN^{\prime}) & = & (-1, \nu)_p^{\binom{v-1}{2}}(-b\nu, v)_p (bv, \nu^{v-1})_p \\
                                            & = & (-1, \nu)_p^{\binom{v-1}{2}}(-b\nu, v)_p (bv, -\delta(G))_p. \end{eqnarray*}
Also,$\epsilon_p(E_0)= \epsilon_p(G)$ by Definition 4.1. Therefore, Lemma 3.8 (b) implies :
\begin{eqnarray*}  1= \epsilon_p(MM^{\prime}) & = &
 \epsilon_p (NN^{\prime}) (-1, \det (E_0))_p\epsilon_p(E_0)\\ & = &
(-1, \nu)_p^{\binom{v-1}{2}}(-b\nu, v)_p (bv, -\delta(G))_p(-1, \delta (G))_p \epsilon_p(G).\end{eqnarray*}
 Therefore,
  \begin{eqnarray*} (-1, \nu)_p^{\binom{v-1}{2}}(\nu,v)_p & = & (-b, v)_p  (bv, -\delta(G))_p(-1, \delta (G))_p \epsilon_p(G)\\
  & = & (b,-v)_p(-bv, \delta (G))_p \epsilon_p (G). \end{eqnarray*}
Since $v=f+1, b=f+g+1$, this proves the second part.$\Box$
\newline

\textbf{Remark.} More generally, in the proof of Theorem 4.2, we might have argued with the most general symmetric non-singular matrix $M$ over $\mathbb{Q}$ containing the incidence matrix $N$ as a sub-matrix. That is, we could take
$$M=  \left( \begin{array}{c} N\\ X \end{array} \right ), $$
where $X$ is any $g \times b$ matrix over $\mathbb{Q}$ for which $M$ is non-singular. Then, as in the above proof, we get $MM^{\prime}/NN^{\prime} = XEX^{\prime}$, where $E$ is as in Definition 4.1. Since $M$ (and
 hence $MM^{\prime}$) is to be non-singular, it follows from Lemma 3.8 that $ XEX^{\prime}$ must be non-singular. Let $X_1,X_2$ be two such choices for $X$, and let $M_1,M_2$ be the corresponding choices for $M$. Since $ (EX_i^{\prime})^{\prime} (EX_i^{\prime})=X_i(E^{\prime}E)X_i^{\prime}= X_i EX_i^{\prime}$ is non-singular, it follows that, for $i=1,2$,  $E_i := EX_i^{\prime}$ is a $b \times g$ matrix of rank $g = \text{rank}(E)$. Therefore, ${\cal C}(E_1) = {\cal C}( E)= {\cal C}(E_2)$. Hence, Lemma 3.9 implies that $X_1EX_1^{\prime} = E_1^{\prime}E_1 \sim E_2^{\prime}E_2 = X_2EX_2^{\prime}$. That is, $M_1M_1^{\prime}/NN^{\prime} \sim M_2M_2^{\prime}/NN^{\prime}$. Hence, by Theorem 3.5, $\det(M_1M_1^{\prime}/NN^{\prime}) \equiv \det (M_2M_2^{\prime}/NN^{\prime})
 \pmod{\mathbb{Q}^{\Box}}$, and $ \epsilon_p(M_1M_1^{\prime}/NN^{\prime})= \epsilon_p(M_2M_2^{\prime}/NN^{\prime})$. Thus $\det (MM^{\prime}/NN^{\prime}) \pmod{\mathbb{Q}^{\Box}}$ and  $ \epsilon_p(MM^{\prime}/NN^{\prime})$ are independent of the choice of $M$. Hence, nothing is to be gained by the apparently most general choice of $M$, and we may as well take the simplest choice, as we have done in the proof.
\newline

\textbf{Application 1 : Strongly resolvable 2-designs.}  A strongly resolvable 2-design may be defined as a quasi-symmetric 2-design with block graph  $K_{m \times n}$, where $m \geq 2, n \geq 2$. The minimal idempotent of rank $g$ is given in this case by : $E=\frac{1}{n} J_n \otimes (I_m- \frac{1}{m} J_m)$. Therefore we may choose the matrix $E_0$ of Definition 3.2 to be $E_0 = \frac{1}{n}(I_{m-1} - \frac{1}{m}J_{m-1})$. Hence, $\det(E_0) = 1/(mn^{m-1})$ and (by Lemma 2.6) $\epsilon_p (E_0) = (-1,n)_p^{\binom{m-1}{2}} (m,n)_p^m (-1,m)_p$.  So we see :
$$ \delta ( K_{m \times n}) = mn^{m-1} \pmod{\mathbb{Q}^{\Box}}, \, \epsilon_p( K_{m \times n})=  (-1,n)_p^{\binom{m-1}{2}} (m,n)_p^m (-1,m)_p. $$
\newline

The following corollary is essentially Theorem 4.1 of \cite{SSS3} (applied to the dual design).

\textbf{Corollary 4.3} : Let $m \geq 2, \, n \geq 2$, and let $\mu$ be positive. Then, for the existence of a strongly resolvable 2-design of defect $\mu$ with block graph $K_{m \times n}$, the following conditions are necessary :
\newline

\begin{flushleft}

(a1) If $m$ is even, then $mn-m+1$ is a perfect square. \\
(a2) If $m$  is odd and $n$ is even, then $(mn-m+1)\mu$ is a perfect square. \\
(a3) If $m$ and $n$ are both odd, then $n(mn-m+1)$ is a perfect square. \\
(b1) If $m \equiv 2 \pmod{4}$ then $n$ is a sum of two squares. \\
(b2) If $ m \equiv 3 \pmod{4}$, $ n \equiv 2 \pmod{4}$, then $ n$ is a sum of two squares. \\
(b3) If $ m \equiv 3 \pmod{4}, n \equiv 0 \pmod{4}$, then $ mn-m+1$ is a sum of two squares.\\
(b4) If $ m \equiv 1 \pmod{4}, n \equiv 2 \pmod{4}$, then $ n(mn-m+1)$ is a sum of two squares. \\
(b5) If $ m, n$ are both odd, and $ m \not \equiv n \pmod{4}$, then the equation $\mu x^2 + (-1)^{\frac{n-1}{2}} ny^2 = z^2$ has
      a non-trivial solution in integers. \\
(b6) If $m \equiv n \equiv 1 \pmod{4}$ then the equation $ nx^2 - \mu y^2 = z^2 $ has a non-trivial solution in integers. \\
(b7) If $ m \equiv n \equiv 3 \pmod{4}$, then, for all prime numbers p,
  $$(-\mu, -n)_p = \begin{cases} -1 \;\text{if}\; \: p=2, \\ 1 \; \text{if} \; p \neq 2. \end{cases}$$
\end{flushleft}

\textbf{Proof} : Using the parameters of  $K_{m \times n}$ as given in the proof of Corollary 2.4, and the new invariants given above, and noting that,  by Theorem 2.1, the order $\nu$
of the design is given by $\nu = n \mu$, the conclusion of Theorem 4.2 reduces to
$$(a)\;  \mu^{m(n-1)}   \equiv  (mn-m+1)n^{mn} \pmod{\mathbb{Q}^{\Box}}, \;\mbox{\text{and}} $$

$$(b) \; (-1, \mu)_p^{\binom{mn-m}{2}} (\mu , mn-m+1)_p  =  (-1,n)_p^{\binom{mn-m}{2}+ \binom{m-1}{2} -1}(n,mn-m+1)_p^{m-1},$$
for all primes $p$.
Clearly, (a) is equivalent to (a1),(a2) and (a3). When $m$ is even, $mn-m+1$ is an odd square, and hence $m(n-1) \equiv 0 \pmod{4}$.
Therefore, in this case, (b) simplifies to $(-1,n)_p^{m/2} = 1$. But Corollary 3.2 shows that this is just (b1). If $m$ is odd and $n$ is even, we have $\mu \equiv mn-m+1 \pmod{\mathbb{Q}^{\Box}}$.  In this case, (b) simplifies to $(-1, mn-m+1)_p^{ (mn-m+1)/2} = (-1,n)_p^{n/2}$.
By Corollary 3.2, this is equivalent to (b2), (b3) and (b4). If $m$ and $n$ are both odd, then $mn-m+1 \equiv n \pmod{\mathbb{Q}^{\Box}}$. In this case, (b) simplifies to $ ((-1)^{(n-1)/2} n, \mu)_p = ((-1,n)_p^{(m+n)/2}$. In view of Theorem 3.1, this amounts to (b5) and (b6) except in the case $m \equiv n \equiv 3 \pmod{4}$. When  $m \equiv n \equiv 3 \pmod{4}$, we get $(-n, \mu)_p(-1,n)_p =1$, i.e., $(-\mu, -n)_p = (-1,-1)_p$
But, by the formula (H5) of Section 3, $(-1,-1)_p =1$ except when $p=2$; $(-1,-1)_2= -1$. This proves (b7). $\Box$.
\newline

Notice that six of the conclusions of Corollary 4.3 rule out certain complete multipartite graphs as possible block graphs. Only (a2), (b5), (b6) and (b7) involve the defect of the design. The conclusion of (b7) can't be rephrased as the solvability of a diophantine equation. By Corollary 2.4,  any Strongly resolvable 2-design of defect $\mu =1$ has $m=n+1$, and is either an affine plane of order $n$ or its complement. In this case, Corollary 4.3 says that if $n \equiv 1\, \mbox{\text{or}} \, 2 \pmod{4}$ then $n$ must be a sum of two squares. Since affine planes of order $n$ are co-extensive with projective planes (symmetric 2-designs with $\lambda =1$) of order $n$, this is just the Bruck-Ryser Theorem of \cite{Bruck}, i.e., the case $\lambda =1$ of Theorem 1.2. It is also easy to see that Part (a) of Corollary 4.3 is equivalent to the corresponding results of Beker in \cite{Beker}.

 Let us say that a quadruple $(\alpha,l,l^*,t)$ of non-negative integers is \textbf{admissible} if $\alpha >0, \alpha (\alpha -1)=ll^*$ and $\alpha$ divides $(l+l^*)^2t$. Recall that Corollary 2.4 gives an explicit bijection between admissible quadruples and the feasible parameters of strongly resolvable $2$-designs. Note that, for given values of $\alpha, \ell, \ell^*$ with $\alpha >0,\ell \ell^* = \alpha (\alpha -1)$, the admissible values of  $t$ vary over the non-negative multiples of $\frac{\alpha}{d}$, where $d$ is the greatest common divisor of $\alpha$ and $(l+l^*)^2$. Corollary 4.4 rules out infinitely many of these  values of $t$, but infinitely many other values survive.

An \textbf{affine resolvable 2-design} is by definition, a strongly resolvable $2$-design with smaller intersection number $\lambda_1=0$. It is immediate from Corollary 2.4 that these have the parameters corresponding to the admissible quadruples $(\alpha,l,l^*,t)=(1,0,n-2,t)$, where $n\geq 2, t \geq 0$. Thus, they have $m=n^2t+n+1, \mu= (n-1)t+1$. Affine resolvable $2$-designs with these parameters have been denoted by $AD(n,t)$ in the literature.

The only known examples of strongly resolvable 2-designs are those obtained in the following construction of Shrikhande and Raghavarao \cite{SSS4} :
\newline

\textbf{Theorem 4.4}:

Suppose there is an affine resolvable 2-design $D_1$ of order $\nu_1$ with block graph $K_{m \times n}$ and a symmetric 2-design $D_2$ of order $\nu_2$ on $n$ points. Then there is a strongly resolvable 2-design $D_1[D_2]$ of order $\nu_1 \nu_2$ with block graph $K_{m \times n}$.
\newline

\textbf{Proof} : For $1 \leq i \leq n $, $ 1 \leq j \leq m $, let $B_{ij}$ be the $i$th block of $D_1$ in its $j$th parallel class (in some order). Without loss of generality, we may assume that the point set of $D_2$ is $\{ 1,2, \ldots,n \}$. Let $D_1[D_2]$ be the incidence system whose blocks are the sets $\cup_{i \in C} B_{ij}$, where $1 \leq j \leq m$, and $C$ varies over the blocks of $D_2$. It is easy to verify that this has the required properties. $\Box$
\newline

Notice that in this construction, $\overline{D_1[D_2]} = D_1[\overline{D_2}]$, so that the class of designs constructed here is closed under complementation. Moreover, when $D_2$ is the $2-(n,1,0)$ design, we get $D_1[D_2]=D_1$, so that this class contains the affine resolvable 2-designs as degenerate cases. Also note that if ($\alpha$ divides $t$ and) there is an affine resolvable $2$-design $D_1=
AD(l+l^* + 2\alpha, \frac{t}{\alpha})$ and a symmetric $2$-design $D_2$ of order $\alpha $ and  balance $\ell $, then the design $D_1[D_2]$ has the  parameters  (given by Corollary 2.4) corresponding to the admissible quadruple $(\alpha,l,l^*, t)$.

When $q$ is a prime power, the design of points versus hyper-planes in the $d$-dimensional affine space $EG(d,q)$ ($d \geq 2$) over the field of order $q$ is an affine resolvable design $AD(n,t)$ with $n=q,\, t =\frac{q^{d-2}-1}{q-1}$. Another series of affine resolvable designs are the $AD(2,t)$ ($t \geq 1$). It is easy to see that these are precisely the $3-(4t+4, 2t +2, t)$ designs (known as the Hadamard $3$-designs). These are co-extensive with Hadamard matrices of order $4t + 4 $, and are expected to exist for all values of $t$.

The only known affine resolvable 2-designs are (a) the Hadamard 3-designs, (b) the affine spaces over finite fields, and (c) other designs having the same parameters as those in (b) and derived from the designs in (b) by algebraic perturbations. Since all known strongly resolvable 2-designs are obtained from these designs via Theorem 4.4, they all  have $n$ prime power and $m \equiv n+1 \pmod{n^2}$. Nobody has bothered to conjecture that $n$ must be a prime power for strongly resolvable 2-designs, since this would include the famous prime power conjecture for projective planes as a very special case. But the following conjecture may be more tractable :

   \textbf{Conjecture}: For the existence of a quasi-symmetric 2-design with block graph $K_{m\times n}$, we must have $m \equiv n+1 \pmod{n^2}$. In other words, we conjecture that the parameters of any strongly resolvable $2$-design corresponding to the quadruple $(\alpha, l, l^*, t)$ can exist only if $\alpha$ divides $t$.

   The smallest feasible parameters of strongly resolvable $2$-designs failing this conjecture correspond to the admissible quadruple $(\alpha,l,l^*,t)= (4,2,6,1)$. These parameters are $n=16, m= 81, \mu =19$,$ b= 1296, v=1216$, $r=486, k=456, \lambda = 182$, and $\lambda_1 = 152, \lambda_2 = 171$. Is there a design with these parameters?
\newline

 \textbf{Application 2: Co-triangular 2-designs.} These are the quasi-symmetric 2-designs with block graph $T_n^*$ ($n \geq 5$).
 The minimal idempotent of rank $g$ in the Bose-Mesner algebra of $T_n^*$ is :

   $$E= \frac{2}{n-2} I_{\binom{n}{2}} - \frac{4}{n(n-2)} J_{\binom{n}{2}} +\frac{1}{n-2} L(K_n), $$
   where $L(K_n)$ denotes the adjacency matrix of $l(K_n)= T_n$. Therefore, the matrix $E_0$ of Definition 4.1 is any non-singular matrix of the form

   $$E_0 =  \frac{2}{n-2} I_{n-1} - \frac{4}{n(n-2)} J_{n-1} +\frac{1}{n-2} L(H),$$
   where $L(H)$ is the adjacency matrix of the line graph $l(H)$ of a graph $H$ with $n$ vertices and $n-1$ edges. For simplicity, we may choose $H= K_{1,n-1}$ (the complete bipartite graph on $1+(n-1)$ vertices), so that $L(H) = J_{n-1}-I_{n-1}$. With this choice, we have :

   $$E_0 = \frac{1}{n-2} (I_{n-1} + \frac{n-4}{n}J_{n-1}). $$
   Thus we get $\det(E_0) =\frac{1}{n(n-2)^{n-3}} $, and (by Lemma 3.6) $\epsilon_p(E_0)= (-1, n-2)_p^{\binom{n-1}{2}} (n-2,n)_p^n (-1,n)_p $. Therefore, by Definition 4.1, we get ;

   $$ \delta (T_n^*)\equiv n(n-2)^{n-1}\pmod{\mathbb{Q}^{\Box}}, \; \epsilon_p(T_n^*)= (-1, n-2)_p^{\binom{n-1}{2}} (n-2,n)_p^n (-1,n)_p. $$

By Corollary 2.5, for any given $\mu \geq 1$, the feasible pairs $(T_n^*, \mu)$ are in bijection with the pairs $(\ell, \ell^*)$ of non-negative integers such that $\ell \ell^* = 4 \mu (\mu -1)$. The correspondence is given by $n= 4 \mu +1 + \ell + \ell^*$ . One may compute using Theorem 2.1 that the parameters corresponding to the pair  $(\ell, \ell^*)$ are as follows :
$$ n= 4 \mu +1 + \ell + \ell^*, \, b=\binom{n}{2}, \, v = \binom{n-1}{2}, $$
$$ r= \frac{n}{2}(\ell + 2 \mu ), \, k= (\frac{n}{2}-1)(\ell + 2 \mu ),\, \lambda = \frac{n}{2}\ell + 2 \mu ,$$
$$ \lambda_1 = (\frac{n}{2}-1)\ell + \mu, \, \lambda_2 = (\frac{n}{2}-1)\ell + 2 \mu.$$

The next corollary is essentially Theorem 5.1 of \cite{SSS3} (applied to the dual design).

\textbf{Corollary 4.5}:

Let $\mu \geq 1$, $n \geq 5$ be integers. For the existence of a co-triangular 2-design of defect $\mu$ with block graph $T_n^*$, the following conditions are necessary :

(a1) If $n \equiv 1\pmod{4}$ then $\mu$ is a perfect square,

(a2) If $n \equiv 2\pmod{4}$ then $(n-2)\mu$ is a perfect square,

(a3) If $n \equiv 3\pmod{4}$ then $n-2$ is a perfect square,

(b1) If $n \equiv 0\pmod{4}$ then $(\mu, (-1)^{n/4}\binom{n-1}{2})_p = ( n-2, (-1)^{n/4}2)_p$

 for all primes $p$,

(b2) If $n \equiv 1\pmod{4}$ then the equation $(n-2)x^2 + (-1)^{(n-1)/4}2y^2 =z^2$ has a non-trivial solution in integers,

(b3) If $n \equiv 2\pmod{4}$ then $n-1$ is a sum of two squares,

(b4) If $n \equiv 3\pmod{4}$ then the equation $\mu x^2 + (-1)^{(n-3)/4} \binom{n-1}{2}y^2 = z^2$ has a non-trivial solution in integers.

\textbf{Proof}

 If $\nu$ is the order of such a design, then Theorem 2.1 gives $\nu = (n-2)\mu$. Therefore, using the spectral parameters  of $T_n^*$ given in the proof of Corollary 2.5, and the new invariants displayed above, we see that, in this case, the conclusions of Theorem 4.2 become :

$$ (a) \; \mu^{\binom{n+1}{2}} \equiv (n-2)^{\binom{n}{2}}\pmod{\mathbb{Q}^{\Box}},$$

$$ (b) \; (-1, \mu)_p^{\binom{\binom{n-1}{2}}{2}} \left(\mu, \binom{n-1}{2}\right)_p = (-1, n-2)_p^{\binom{\binom{n-1}{2}+1}{2} + n-1} (2,n-2)_p  $$

for all primes $p$. (To verify that Part (b) of Theorem 4.2 reduces to (b) above, we need the following formula : $(n, n-2)_p (-2, n)_p = (2, n-2)_p$. Proof of this formula : $(n,n-2)_p = (2,2)_p (2, n/2)_p(2, n/2-1)_p = (-2, n/2)_p (2, n/2 -1)_p = (-2,n)_p(2,n-2)_p$ since $(2,2)_p = (2,-2)_p(1-(-1), -1)_p=1$.)

 Now, (a) is clearly equivalent to (a1), (a2) and  (a3). When $n \equiv 0 \pmod{4}$, we have $ \binom{\binom{n-1}{2}}{2}\equiv \frac{n}{4} \pmod{2}$ and $\binom{\binom{n-1}{2}+1}{2} \equiv \frac{n}{4}-1 \pmod{2}$. Hence, (b) reduces to (b1) in this case.

 When $n \equiv 1 \pmod{4}$, we have $\mu \equiv 1 \pmod{\mathbb{Q}^{\Box}}$ and  $\binom{\binom{n-1}{2}+1}{2} \equiv \frac{n-1}{4}\pmod{2}$,
 so that (b) reduces to $(n-2, (-1)^{\frac{n-1}{4}}2)_p =1$ for all primes $p$. When $n \equiv 2 \pmod{4}$, we have $\mu \equiv n-2  \pmod{\mathbb{Q}^{\Box}}$ and $ \binom{\binom{n-1}{2}}{2} \equiv \binom{\binom{n-1}{2}+1}{2} \pmod{2}$. Hence (b) reduces to $(-1, n-1)_p =1$ for all primes $p$ When $n \equiv 3 \pmod{4}$, we have $n-2 \equiv 1  \pmod{\mathbb{Q}^{\Box}}$ and  $ \binom{\binom{n-1}{2}}{2} \equiv \frac{n-3}{4} \pmod{2}$. Therefore, in this case, (b) reduces to $(\mu, (-1)^{\frac{n-3}{4}} \binom{n-1}{2})_p =1$ for all primes $p$. In view of Lemma 3.1 and Corollary 3.2, these observations complete the proof. $\Box$

By Corollary 2.5 (1), the co-triangular designs of defect $1$ with block graph $T_n^*$  are co-extensive with biplanes of order $n-2$. In this case, Corollary 4.5 reduces to the case $\lambda = 2$ of Theorems 1.1, 1.2. (Notice that when $n-2$ is a square, $n-1$ is trivially a sum of two squares.) We can verify that Corollary 4.5(a) is equivalent to the result in \cite{Haemers}. No co-triangular designs of defect $\mu >1$ are known.
\newline

\textbf{Application 3 : The symplectic designs.} These are the quasi-symmetric 2-designs with block graph $Sp(2d,2)$ ($d \geq 3$), the non-orthogonality graph of a $(2d)$-dimensional symplectic space over the field of order $2$. In view of Corollary 2.7 and the remark following it, we ignore the case $q>2$.
Following the discussion in \cite{Cameron}, p. 75, two  series of quasi- symmetric 2-designs  may be obtained by dualizing the contraction at a point of the remarkable symmetric 2-designs (or their complements) described in Example 5.17 of \cite{Cameron}, p. 76. The designs in the first series have block graph $Sp(2d,2)$, and may be described directly as follows (the second series, with block graph $Sp(2d,2)^*$, has a similar description with the elliptic quadrics replaced by hyperbolic quadrics).

\textbf{A construction of quasi-symmetric 2-designs} with block graph $Sp(2d,2)$. Fix $d \geq 3$. The blocks of the design are the points of the projective space $PG(2d-1,2)$. The points are all the elliptic quadrics (irreducible quadrics of Witt index $d-1$ ) in $PG(2d-1,2)$ whose defining quadratic forms polarise to a given non-degenerate symplectic bilinear form on the underlying vector space. The incidence is reverse containment. These designs have the following parameters :

$$ b= 2^{2d}-1, \, v= 2^{d-1}(2^d -1), \, r= (2^d +1)(2^{d-1}-1), \, k = 2^{d-1}(2^{d-1}-1), $$
$$ \lambda = 2^{2d-2} - 2^{d-1} -1, \, \lambda_1 = 2^{d-1}(2^{d-2}-1), \, \lambda_2 = 2^{d-2}(2^{d-1} -1) .$$

 Applying Theorem 4.2 to these designs, we deduce that
 $$ \delta (Sp(2d,2))= 2^{d-1}(2^d+1)\pmod{\mathbb{Q}^{\Box}},\; \epsilon_p( Sp(2d,2))= (2, 2^{2d}-1)_p^d. $$

(To deduce this compact formula for $\epsilon_p(\cdot)$, we need to observe that, for any rational number $x \neq \pm 1$, $(1-x,1+x)_p =
 (2, 1- x^2)_p $. Proof : we have $1 = (\frac{1}{2}(1-x), \frac{1}{2} (1+x))_p =(1-x,1+x)_p (2, 1-x^2)_p $.)
\newline

\textbf{Corollary 4.6}

For the existence of a quasi-symmetric 2-design of order $\nu$ with block graph $Sp(2d,2)$, $\nu$ must be a perfect square.
\newline

\textbf{Proof} :

The spectral parameters of $ Sp(2d,2)$ are
$$ \rho = 2^{d-1}, \, \sigma = - 2^{d-1}, \, f = (2^d + 1)(2^{d-1}-1), \, g = (2^d - 1)(2^{d-1} +1).$$

Combining these parameters with the formulae for the discriminant and the p-adic invariants of $ Sp(2d,2)$ given above, we see that, in this case, the first part of Theorem 4.2 reduces to $\nu \equiv 1 \pmod{\mathbb{Q}^{\Box}}$. (In this case, the second part of Theorem 4.2 follows from the first part.) $\Box $
\newline

\textbf{Application 4 : The multi-Steiner designs.} These are the quasi-symmetric 2-designs with block graph $S_n(m)$ ( $2 \leq n <m, n|m(m-1)$).
Applying Theorem 4.2 to the Steiner 2-design ($\mu =1$) whose block graph is $S_n(m)$, we may easily deduce the formulae for the discriminant and p-adic invariant of the Steiner graphs :

$$ \delta(S_n(m))= (m-1)^{m(n-1)} \pmod{\mathbb{Q}^{\Box}},$$
$$ \epsilon_p(S_n(m))= (-1, m-1)_p^{\binom{mn-m}{2}-1}(-mn, m-1)_p^{mn-m}(mn(m-1), -(mn-m+1))_p. $$
\newline

\textbf{Corollary 4.7} : Let $m> n \geq 2$ be integers such that $n$ divides $ m(m-1)$. Let $\mu$ be a positive integer. Then, for the existence of a multi-Steiner 2-design of defect $\mu$ and block graph $ S_n(m)$, the following conditions are necessary :
\newline

\begin{flushleft}
(a) If $m$ is odd and $n$ is even, then $\mu$ must be a perfect square,

(b) The diophantine equation $\mu x^2 + (-1)^{\binom{mn-m}{2}}(mn-m+1)y^2 = z^2$ must have a non-trivial solution in integers $x,y,z.$
\end{flushleft}

\textbf{Proof} :

In view of the spectral parameters (as given in the proof of Corollary 2.8), and the discriminant and p-adic invariants of $S_n(m)$ given above, the conditions in Theorem 4.2 reduce in this case to :
$$ (a)\; \mu^{m(n-1)} \equiv 1 \pmod{\mathbb{Q}^{\Box}}, \newline (b) \, (\mu,  (-1)^{\binom{mn-m}{2}}(mn-m+1))_p=1. $$

By Lemma 3.1, this is just the conclusion of this Corollary. $\Box$
\newline

\textbf{Remark} : Let $ d \geq 3$, $\Gamma_d = Sp(2d,2)$ and let $\Lambda_d$ be the block graph of a $2-(2^{d-1}(2^d -1), 2^{d-1}, 1)$ design (that is, $\Lambda_d = S_n(m)$ where $n=2^{d-1}$, $m = 2^d +1$). (Such designs, and hence graphs, exist for all values of $d$. An example is the incidence system whose points are the lines of  $PG(2, 2^d)$ disjoint from a given hyper-oval in the projective plane, blocks are the points of the plane outside the hyper-oval, and incidence is reverse containment.) Observe that $\Gamma_d$ and $\Lambda_d$ have the same parameters, but
$$ \delta (\Gamma_d) =  2^{d-1}(2^d+1)\pmod{\mathbb{Q}^{\Box}},\; \delta (\Lambda_d) =  2(2^d+1)\pmod{\mathbb{Q}^{\Box}},$$
$$ \epsilon_p(\Gamma_d)= (2, 2^{2d}-1)_p^d, \; \epsilon_p(\Lambda_d)=1.$$

Thus, when $d$ is an odd number, $\delta(\Gamma_d) \neq \delta(\Lambda_d)$. Also, if $d$ is odd and $p \equiv \pm 3 \pmod{8}$  is a prime number dividing the square-free part of $2^{2d}-1$, then $\epsilon_p(\Gamma_d) \neq \epsilon_p(\Lambda_d)$. For example, whenever $d \equiv \pm 1 \pmod{6}$, we get $\epsilon_3(\Gamma) \neq \epsilon_3(\Lambda)$. So the usual parameters of an s.r.g. do not determine its new invariants.

Since isomorphic graphs clearly have the same invariants, and since  $\delta(\Gamma_d) \neq \delta(\Lambda_d)$  for odd $d$, it follows that $\Gamma_d$ and $\Lambda_d$ are non-isomorphic for all odd $d$. Since this holds for any Steiner graph $\Lambda_d$ with the parameters of $\Gamma_d$, we have :
\newline

\textbf{Corollary 4.8} :

The symplectic graphs $Sp(2d, 2)$ are not geometrizable for odd numbers $d \geq 3$. That is, these graphs are not block graphs of Steiner $2$-designs.

However, $Sp(4,2)=T_6$ is geometrizable : it is the block graph of the $2-(6,2,1)$ design.\\
\newline
 \textbf{Question :} For what (even) values of $d$ is the graph $Sp(2d,2)$ geometrizable ?
\newline

The Steiner 2-designs and their complements are the only multi-Steiner designs with $\mu =1$. The only known examples of non-trivial  multi-Steiner designs are the designs $PG_{d-2}(d, q)$ ($d \geq 4$, $q$ prime power) whose points and blocks are the points and $(d-2)$-dimensional flats in the $d$-dimensional projective space $PG(d,q)$ over the field of order $q$. The block graph of this design is isomorphic to that of the Steiner 2-design  $PG_1(d,q)$ (the design of points versus lines in $PG(d,q)$.) Any duality of the projective space induces an isomorphism between these two graphs.

In lieu of a convenient description of the feasible parameters of multi-Steiner 2-designs, we present below a table of
small parameters. Only the smaller of a pair of complementary parameters is given. An entry "no" in the "exists ?" column means that it is ruled out by Corollary 4.7. (Note that this corollary does not rule out designs with a pseudo-Steiner block graph! For instance, a design with the parameter of item number $10$ in Table 1 does exist, but with block graph $Sp(6,2)$.) An "yes" entry here means that the parameters are in the series given in the previous paragraph (the only construction we know!).

It is conceivable that item number 12 of Table 1 exists with the block graph of $EG_1(3,4)$. But our preliminary investigation makes it look unlikely. A much more promising candidate is the block graph of the classical unital with automorphism group $U(3,5)$ as the block graph for item 24 of this table.

\begin{table}
\centering
\caption{Small feasible parameters of q.s. designs with block graph $S_n(m)$}

\begin{tabular}{ccccccccc}
\hline
Number & $n$  & $m$ & $v$ & $k$ & $\lambda$ & $\lambda_1$ & $\lambda_2$ & exists? \\
\hline
1       & $3$ & $10$ & $21$ & $9$ & $12$  & $3$ & $5$   & no \\
2        & $3$ & $15$ &$31$  & $7$ & $7$   & $1$ & $3$ & yes \\
3        & $3$ & $16$  & $33$ & $15$ & $35$ & $6$ & $9$  & ? \\
4        & $3$ & $19$  & $39$ & $12$ & $22$  & $3$ & $6$  & ? \\
5        & $3$ & $22$ & $45$ & $21$ & $70$   & $9$ & $13$  & ? \\
6        & $3$ & $27$ & $55$ & $16$ & $40$ & $4$ & $8$ &  ? \\
7        & $3$ & $31$ & $63$ & $15$ & $35$ & $3$ & $7$ & yes \\
8        & $3$ & $36$ & $73$ & $10$ & $15$ & $1$ & $4$ & ? \\
9        & $3$ & $66$ & $133$ & $13$ & $26$ & $1$ & $5$ & ? \\
10       & $4$ & $9$ & $28$    & $12$ & $11$ & $4$ & $6$ & No \\
11       & $4$ & $17$ & $ 52$  & $16$ & $20$ & $4$ & $7$ & No \\
12       & $4$ & $21$ & $64$  & $24$ & $46$ & $8$ & $12$ & ? \\
13        & $4$ & $40$ & $121$ & $13$ & $13$ & $1$ & $4$ & yes \\
14        & $5$ & $16$ & $65$  & $20$ & $19$ & $4$ & $7$ & No \\
15        & $5$ & $26$ & $105$ & $25$ & $30$ & $5$ & $9$ & ? \\
16         & $5$ & $45$ & $181$ & $16$ & $12$ & $1$ & $4$ & ? \\
17         & $5$ & $85$ & $341$ & $21$ & $21$ & $1$ & $5$ & ? \\
18         & $6$  & $9$ & $46$  & $16$ & $8$ & $4$ & $6$ & No \\
19         & $6$  & $10$ & $51$ & $15$ & $7$ & $3$ & $5$ & No \\
20         & $6$  & $13$ & $66$  & $30$ & $29$ & $12$ & $15$ & No \\
21         & $6$  & $18$ & $91$ & $40$ & $52$ & $16$ & $20$ & ? \\
22         & $6$  & $19$ & $96$ & $36$ & $42$ & $12$ & $16$ & ? \\
23         & $6$  & $22$ & $111$ & $21$ & $14$ & $3$ & $6$ & ? \\
24         & $6$  & $25$ & $126$ & $30$ & $29$ & $6$ & $10$ & ? \\
25         & $6$   & $96$ & $481$  & $25$  & $20$ & $1$ & $5$ & ? \\
\hline
\end{tabular}
\end{table}

\end{document}